\documentclass[preprint,12pt]{elsarticle}

\usepackage{amssymb}
\usepackage{amsmath}
\usepackage{float}
\usepackage{graphicx}
\usepackage{booktabs}
\usepackage{subfigure}
\usepackage{graphbox}
\usepackage{hyperref}
\usepackage{siunitx}
\usepackage{upgreek}
\usepackage{amsmath}
\usepackage[nameinlink]{cleveref}
\crefname{figure}{Fig.}{Figs.} % capitalize "F", no period
\crefname{equation}{Eq.}{Eqs.} % capitalize "E", no period
\crefname{table}{Table}{Tables} % capitalize "T", no period
\crefname{section}{Section}{Sections} % capitalize "T", no period
\crefrangelabelformat{equation}{(#3#1#4--#5#2#6)}
\usepackage{bm}
\usepackage[section]{placeins}
\usepackage{tikz}
\usetikzlibrary{positioning, calc, matrix, arrows.meta}

\journal{CMAME}

\begin{document}

\begin{frontmatter}

%% Title, authors and addresses

%% use the tnoteref command within \title for footnotes;
%% use the tnotetext command for theassociated footnote;
%% use the fnref command within \author or \affiliation for footnotes;
%% use the fntext command for theassociated footnote;
%% use the corref command within \author for corresponding author footnotes;
%% use the cortext command for theassociated footnote;
%% use the ead command for the email address,
%% and the form \ead[url] for the home page:
%% \title{Title\tnoteref{label1}}
%% \tnotetext[label1]{}
%% \author{Name\corref{cor1}\fnref{label2}}
%% \ead{email address}
%% \ead[url]{home page}
%% \fntext[label2]{}
%% \cortext[cor1]{}
%% \affiliation{organization={},
%%            addressline={}, 
%%            city={},
%%            postcode={}, 
%%            state={},
%%            country={}}
%% \fntext[label3]{}

\title{Efficient thermal simulation in metal additive manufacturing via semi-analytical isogeometric analysis}
%% Article title

%% use optional labels to link authors explicitly to addresses:
% Author names
\author[address_ME]{Yang Yang}
\ead{Y.Yang-13@tudelft.nl}
\ead[orcid]{0000-0002-2586-9218}

\author[address_AM]{Ye Ji \corref{cor1}}
\cortext[cor1]{Corresponding author}
\ead{y.ji-1@tudelft.nl}
\ead[orcid]{0000-0002-1173-6457}

\author[address_AM]{Matthias Möller}
\ead{m.moller@tudelft.nl}
\ead[orcid]{0000-0003-0802-945X}

\author[address_ME]{Can Ayas}
\ead{C.Ayas@tudelft.nl}
\ead[orcid]{0000-0002-9472-7424}

% Author affiliation
\affiliation[address_ME]{organization={Computational Design and Mechanics Group, Department of Precision and Microsystems Engineering, Faculty of Mechanical Engineering},
    addressline={Mekelweg 2}, 
    city={Delft},
    % citysep={}, % Uncomment if no comma needed between city and postcode
    postcode={2628 CD}, 
    % state={},
    country={The Netherlands}}

\affiliation[address_AM]{organization={Delft Institute of Applied Mathematics, Delft University of Technology},
    addressline={Mekelweg 4},
    city={Delft},
    % citysep={}, % Uncomment if no comma needed between city and postcode
    postcode={2628 CD}, 
    % state={},
    country={The Netherlands}}

%% Abstract
\begin{abstract}

Thermal modeling of Laser Powder Bed Fusion (LPBF) is challenging due to steep, rapidly moving thermal gradients induced by the laser, which are difficult to resolve accurately with conventional Finite Element Methods (FEM). Highly refined, dynamically adaptive spatial discretization is typically required, leading to prohibitive computational costs. Semi-analytical approaches mitigate this by decomposing the temperature field into an analytical point-source solution and a complementary numerical field that enforces boundary conditions. However, state-of-the-art implementations either necessitate extensive mesh refinement near boundaries or rely on restrictive image-source techniques, limiting their efficiency and applicability to complex geometries. 
This study presents a novel reformulation of the semi-analytical framework using Isogeometric Analysis (IGA). The laser heat input is captured by the analytical point-source solution, while the complementary correction field, which imposes boundary conditions, is solved using a spline-based IGA discretization. The governing heat equation for the correction field is cast in a weak form, discretized with NURBS basis functions, and advanced in time using an implicit $\theta$-scheme. This approach leverages IGA's key advantages: exact geometry representation, higher-order continuity, and superior accuracy per degree of freedom. These features unlock efficient thermal modeling of realistic parts with complex contours. 
Our strategy eliminates the need for scan-wise remeshing and robustly handles intricate geometric features like sharp corners and varying cross-sections. Numerical examples demonstrate that the proposed semi-analytical IGA method delivers accurate temperature predictions and achieves substantial computational efficiency gains compared to standard FEM, establishing it as a powerful new tool for high-fidelity thermal simulation in LPBF.
\end{abstract}

% %%Graphical abstract
% \begin{graphicalabstract}
% %\includegraphics{grabs}
% \end{graphicalabstract}

%%Research highlights
% \begin{highlights}
% \item Research highlight 1
% \item Research highlight 2
% \end{highlights}

%% Keywords
\begin{keyword}
% %% keywords here, in the form: keyword \sep keyword

% %% PACS codes here, in the form: \PACS code \sep code

% %% MSC codes here, in the form: \MSC code \sep code
% %% or \MSC[2008] code \sep code (2000 is the default)
Heat transfer modeling \sep Isogeometric analysis \sep Semi-analytical method \sep Laser powder bed fusion \sep Additive manufacturing \sep Complex geometries
\end{keyword}

\end{frontmatter}

%% Add \usepackage{lineno} before \begin{document} and uncomment 
%% following line to enable line numbers
%% \linenumbers

%% main text
%%

%% Use \section commands to start a section
\section{Introduction}
\label{sec:introduction}

Laser Powder Bed Fusion (LPBF) is a metal additive manufacturing (AM) technology that offers exceptional form freedom to fabricate geometrically complex parts. It has been widely adopted for manufacturing topology-optimized structures and meta-materials \cite{BAYAT2023101129,LIU2025124941,BIFFI2024100216}. In LPBF, one or multiple laser beams are guided by mirrors and focused onto the metal powder bed. The absorbed laser energy locally melts the powder to form small melt pool(s). As the laser(s) continuously move, the melt pool(s) rapidly cool and solidify, realizing the layer-by-layer growth. During this, powder–liquid–solid phase changes, and several thermal phenomena arise. Heat conduction within the solid, from the localized melt pool toward the build platform acting as a heat sink, accounts for conveying most of the heat input through laser(s). In contrast, heat transfer to the surrounding loose powder dissipates only about 1\% of the heat conducted \cite{gusarov2009model}, while radiation and convection to the surroundings result in only minor heat losses.

Accurate thermal simulations play an important role in predicting the temperature histories of additively manufactured parts. Manufacturing parameters such as laser power, scanning speed, and layer thickness can be optimized based on the temperature histories. Moreover, reliable temperature fields inform models of phase transformations, microstructure evolution, residual stress, and distortion, and ultimately improve part quality.

The Finite Element Method (FEM) remains the most widely used approach for part-scale temperature simulations in LPBF. The temperature fields obtained from these simulations have been used for predicting thermal stress \cite{JIMENEZABARCA2023110151,YANG2021102090}, porosity formation mechanisms \cite{VASTOLA2018817} and melt-pool shape \cite{trejos2022finite}. Comprehensive review of state-of-the-art FEM applications can be found in \cite{SARKAR2024104157}.  
To represent the layer-by-layer growth of a part during the AM process, FEM typically employs element birth–death techniques. Two common variants are the quiet element method, where undeveloped layers remain in the global conductivity matrix but are rendered thermally inactive by drastically reducing their conductivity \cite{olleak2020part}, and the inactive element method, where the elements associated with these layers are excluded entirely until activated during deposition \cite{ROBERTS2009916}. 

Multi scale nature of the LPBF problem imposes challenges on FEM based process simulations. Accurate resolution of steep temperature gradients near a moving laser spot demands extremely fine meshes on the order of the laser spot radius \cite{ROBERTS2009916,FOROOZMEHR2016255}. For example, Roberts et al. \cite{ROBERTS2009916} used \SI{25}{\micro\meter} elements for a \SI{50}{\micro\meter} spot radius, Vanini et al. \cite{VANINI2024104369} refined further to \SI{20}{\micro\meter}, and Zhang et al. \cite{ZHANG2023108839} proposed \SI{5}{\micro\meter} elements for similar conditions. While such fine discretization is required for accuracy, it dramatically increases the total number of elements and thus the computational cost, limiting the practicality of FEM for part-scale LPBF simulations.

Rather than utilizing a uniformly fine mesh across the entire computational domain, adaptive mesh refinement improves computational efficiency while retaining the accuracy of numerical simulations by dynamically adjusting the spatial resolution in response to the changing position of the laser spot. In this approach, high-resolution meshes are employed in the laser–material interaction region, while the mesh is progressively coarsened in areas distant from the heat source. This strategy enables the efficient resolution of steep, localized temperature gradients induced by the moving heat source, which would otherwise necessitate a globally fine mesh and incur substantial computational cost. Within adaptive mesh refinement, two principal refinement strategies are commonly employed: $h$-refinement (reducing the mesh size $h$) and $hp$-refinement (refining $h$ while simultaneously increasing the polynomial order $p$ of the approximation space) \cite{GOUGE2019100771,KOLLMANNSBERGER20181483,LEONOR2024116977,moreira2022multi}. 
Moreover, generally, there are two types of $h$ refinement applied. The more prevalent is the layer-wise (static) remeshing scheme, in which the refinement focuses on the active layer while keeping the underlying previously solidified layers coarser \cite{GOUGE2019100771,LI2019100903,RAUNER2025118769}. The other is the scan-wise (dynamic) remeshing, which refines the mesh locally around the moving laser spot while maintaining a coarser discretization in regions of nearly uniform temperature \cite{OLLEAK202075,OLLEAK2022100051}. 
The layer-wise remeshing strategy primarily focuses on mesh adaptation along the build direction (orthogonal to the layer), ensuring adequate resolution across the successively deposited layers. In contrast, the scan-wise remeshing strategy refines the mesh within the current printing plane, allowing for a more detailed physical representation of the thermal process in the plane. Besides, the scan-wise remeshing strategy can employ finer local meshes, thereby enhancing the accuracy of melt pool size, shape, and cooling rate predictions. However, the scan-wise strategy also introduces greater computational complexity, as it requires much more frequent updates of the mesh and the corresponding thermal conductivity matrix and thermal load vector. Olleak and Xi \cite{OLLEAK202075} optimized the remeshing region length and corresponding remeshing frequency in the scan-wise adaptive remeshing strategy; however, even the most effective configuration achieves only about a 60\% reduction in total computational time relative to the layer-wise remeshing approach. Thus, in general, most of the simulation tools are based on the layer-wise remeshing strategy.

In addition to FEM, isogeometric analysis (IGA) \cite{cottrell2009isogeometric} has emerged as a powerful alternative for thermal simulation in geometrically complex domains \cite{duvigneau2009introduction,zang2021isogeometric,gupta2025transient} and for moving heat sources in LPBF \cite{ji2022curvature}. By employing spline-based basis functions such as NURBS \cite{borden2011isogeometric} or T-splines \cite{bazilevs2010isogeometric,dorfel2010adaptive,scott2011isogeometric,evans2015hierarchical}, IGA preserves the exact CAD geometry, whereas FEM often degrades geometric fidelity during mesh generation. This exact representation, combined with the higher-order continuity of IGA's basis functions, enables accurate temperature predictions with significantly fewer degrees of freedom than conventional FEM \cite{yu2020locally}, making it particularly well-suited for capturing the intricate features of LPBF components \cite{hughes2005isogeometric}. Carraturo et al.~\cite{CARRATURO2019660,CARRATURO2024133} demonstrated that THB-spline refinement can dynamically remesh and track the moving laser spot within an IGA framework—similar in concept to scan-wise remeshing in FEM—with element size ratios ranging from \num{20} to \num{100}. However, despite these adaptive strategies, both FEM and IGA typically require highly refined meshes and frequent updates to resolve the steep, laser-induced thermal gradients, which incurs a significant computational cost. A critical challenge in this context is the substantial time consumed by two core stages: matrix assembly and the solution of the linear system. While IGA offers superior accuracy, its high-order, high-continuity basis functions lead to increased assembly costs per element and produce linear systems with larger bandwidths and poorer conditioning. This often makes the solution phase the dominant computational bottleneck, particularly in large-scale adaptive simulations \cite{schillinger2012isogeometric}. Consequently, extensive research has focused on developing fast assembly techniques \cite{bressan2019sum,pan2020fast} and efficient, specialized solvers \cite{gahalaut2013multigrid,carson2020cost}.

Another approach to thermal modeling in LPBF employs analytical solutions of the heat equation in infinite or semi-infinite domains \cite{YANG2018284,NING2019319,YANG2020100955,STEUBEN2019437,WOLFER2019100898,YANG2025127059}. These methods model the movement of the laser by sequentially activating discrete point heat sources. The closed-form solutions inherently capture the steep temperature gradients near the laser spot, eliminating the need for adaptive mesh refinement and reducing computational complexity. However, since the analytical solution assumes a semi-infinite medium, its accuracy is limited when applied to finite parts, as the physical boundary conditions are not satisfied \cite{YANG2025127059}.

To enforce the correct boundary conditions, semi-analytical methods superpose the analytical field with a complementary numerical field \cite{YANG2018284}. Since the point-source solution captures the steep thermal gradients, the complementary field remains smooth, allowing the numerical mesh to be decoupled from the laser spot size. This permits a significantly coarser discretization without sacrificing accuracy. However, a key limitation arises near the part boundaries: the outgoing heat flux from an analytical source close to a boundary must be counteracted by the numerical field, which in turn demands local mesh refinement dictated by the laser spot scale, reintroducing computational cost.

\begin{figure}[H]
  \centering
  \includegraphics[width=0.9\linewidth]{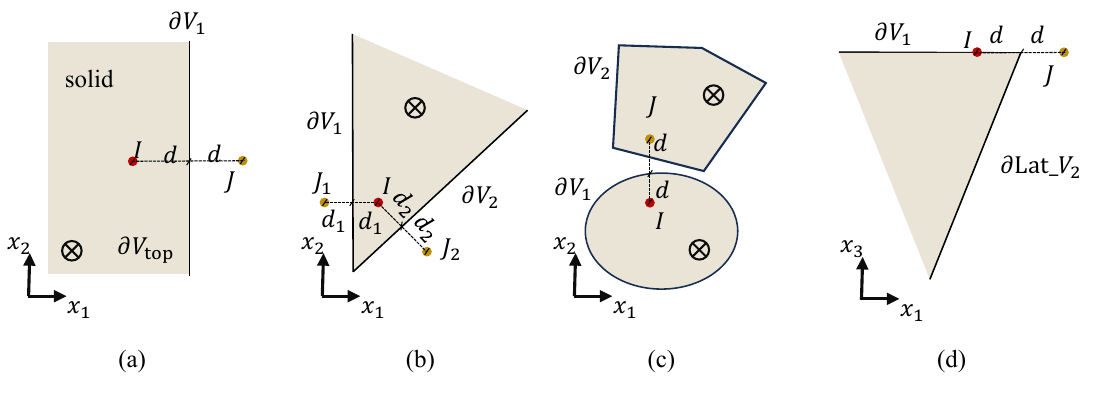}
  \caption{Schematic illustration of the image source: (a) Single straight boundary, for imaging through a simple straight boundary $\partial V_1$, the adiabatic boundary conditions can be easily satisfied by introducing a single image source $J$.
  (b) Multiple connected boundaries, in cases involving multiple connected boundaries, the image source $J_1$ reflected from one boundary $\partial V_1$ may also influence another boundary $\partial V_2$, and similarly, the image source $J_2$ reflected from $\partial V_2$ may affect $\partial V_1$.
  (c) Two adjacent solid boundaries, when imaging between two neighboring solid boundaries, the image source $J$ reflected by one boundary $\partial V_1$ may be placed within the solid domain associated with the other boundary.
  (d) Cross-sectional variation along the building direction ($x_3$), when the cross section varies along the $x_3$ direction, the image source $J$ can only satisfy the adiabatic boundary conditions on the top plane. In this case, the physical source $I$ and its image $J$ are separated by an equal distance $d$ to the boundary only on the top plane.}
\label{fig: figure_boundary_image_problem}
\end{figure}

To circumvent this boundary refinement, our previous work employed analytical image sources across adiabatic boundaries \cite{YANG2018284, YANG2025127059}. For a single straight boundary, the method is straightforward: an image source $J$ is positioned symmetrically to the actual source $I$ across the boundary $\partial V_1$, satisfying the adiabatic condition exactly (\cref{fig: figure_boundary_image_problem}a). For curved boundaries, the implementation becomes more complex \cite{YANG2025127059}.

The image source method, however, exhibits fundamental limitations for the complex geometries typical of LPBF. As illustrated in \cref{fig: figure_boundary_image_problem}b, at sharp corners where multiple boundaries connect, reflecting source $I$ across $\partial V_1$ generates image $J_1$, which itself influences $\partial V_2$. This, in turn, requires a new image $J_2$ from $\partial V_2$, leading to an infinite series of reflections to satisfy all boundary conditions simultaneously. Furthermore, the method fails for other common geometric features: an image source may be incorrectly placed inside an adjacent solid domain, artificially heating it (\cref{fig: figure_boundary_image_problem}c); or for parts with varying cross-sections, the adiabatic condition is satisfied only on the top plane, but not on the subsurface layers (\cref{fig: figure_boundary_image_problem}d). Consequently, while effective for simple geometries, the image source method lacks the generality required for realistic LPBF components.

This paper introduces a novel semi-analytical framework for LPBF thermal simulation that synergistically combines the strengths of the analytical solution and IGA. We reformulate the boundary correction problem by replacing the restrictive image-source technique and standard FEM with an IGA-based solver for the complementary field. This approach delivers two fundamental benefits:
\begin{itemize}
    \item The analytical solution resolves the steep, laser-scale thermal gradients, eliminating the need for a fine, laser-following mesh and frequent dynamic remeshing.
    \item IGA provides a geometrically exact and high-continuity discretization, enabling accurate and efficient enforcement of boundary conditions on complex parts, even with a coarse mesh.
\end{itemize}
These properties collectively enable accurate temperature predictions for realistic parts with complex geometries while requiring substantially fewer degrees of freedom than conventional FEM. Numerical examples demonstrate superior computational efficiency and accuracy compared to existing semi-analytical FEM approaches from the literature.

The remainder of this paper is organized as follows. \Cref{sc: Mathematical Model} presents the reformulation of the semi-analytical method based on discrete point sources, incorporating an IGA-enabled boundary-correction strategy, and further details the IGA discretization and the time-integration procedure employed to compute the numerical complementary field. \Cref{sc: numerical examples} evaluates the efficiency and accuracy of FEM and IGA under varying mesh densities for a point-source, and further demonstrates the robustness of the proposed method for continuous laser scanning and geometrically complex parts. Finally, \Cref{sc:conclusions} summarizes the key findings and contributions.

\section{IGA-based formulation of the semi-analytical method} 
\label{sc: Mathematical Model}

\begin{figure}[H]
  \centering
  \includegraphics[width=0.9\linewidth]{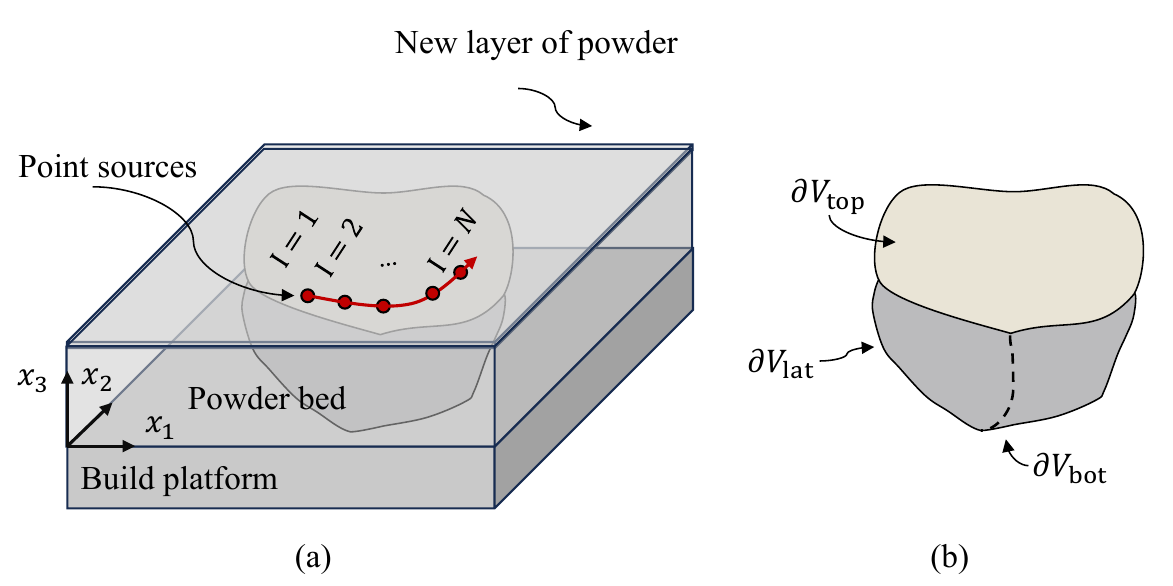}
  \caption{Schematic of the LPBF process: (a) Body $V$ submerged in the powder bed, with bottom $\partial V_{\mathrm{bot}}$ fused to the build platform and a thin powder layer on $\partial V_{\mathrm{top}}$; (b) Boundary decomposition of $V$, showing $\partial V_{\mathrm{bot}}$, $\partial V_{\mathrm{lat}}$, and $\partial V_{\mathrm{top}}$.}
  \label{fig: physical model}
\end{figure}

We consider the LPBF process of a freeform three-dimensional part. To print a new layer, a thin layer of metal powder is recoated over the previously built body $V$ (see \cref{fig: physical model}) with the coordinate origin defined in \cref{fig: physical model}a. The boundary $\partial V$ of $V$ consists of three surfaces: the bottom surface $\partial V_{\mathrm{bot}}$ fused to the build platform, the lateral surface $\partial V_{\mathrm{lat}}$ surrounded by loose powder, and the top surface $\partial V_{\mathrm{top}}$ covered by the newly recoated layer (\cref{fig: physical model}b). When laser scans $\partial V_{\mathrm{top}}$ to fuse the powder, partial absorption of laser energy by the powder drives transient heat conduction within the solid part, governed by  
\begin{equation}\label{eq:heat-conduction}
  \rho c_p \frac{\partial T}{\partial t} = \nabla \cdot (k \nabla T) + Q, \quad \text{in } V,
\end{equation}
where $T(\mathbf{x},t)$ is temperature, $\rho$ is mass density, $c_p$ is specific heat capacity, $k$ is thermal conductivity, and $Q$ is the volumetric heat source.  
Neglecting the temperature dependence of $\rho$, $c_p$, and $k$, the equation becomes  
\begin{equation}\label{eq:linear-heat}
  \frac{\partial T}{\partial t}
  = \alpha \nabla^{2}T + \frac{Q}{\rho c_p},
\end{equation}
where $\alpha = \frac{k}{\rho c_p}$ is the linear form with thermal diffusivity.

The thermal source $Q$ is represented as a superposition of $N$ instantaneous point sources obtained by discretizing the continuous laser scan (see \cref{fig: physical model}a).  
The $I$-th source is activated at time $t^{(I)}$ and the next at $t^{(I+1)} = t^{(I)} + \Delta t$; thus the arc-length spacing between successive sources along the scan path is $v\,\Delta t$, where $v$ denotes the laser scanning speed. Guided by our previous convergence study \cite{YANG2018284,YANG2025127059}, we set $\Delta t=\SI{1e-5}{s}$ for all simulations. After completing the scan of the current layer, a fresh powder layer is recoated.

The temperature field $T(\mathbf{x},t)$ within the domain $V$ under the thermal load $Q$ in \cref{eq:linear-heat} is obtained by solving a boundary value problem with prescribed conditions on $\partial V$.  
The initial condition $t=0$ is  
\begin{equation}
    T(\mathbf{x},0)=T_c,
\end{equation}  
where $T_c$ is the constant temperature of the build platform.  

Because the thermal conductivity of the surrounding powder is only about $1\%$ of that of the solid material \cite{gusarov2009model}, the lateral surface $\partial V_{\mathrm{lat}}$ (see \cref{fig: physical model}) is treated as adiabatic \cite{YANG2020100955,STEUBEN2019437}.  
Heat loss by convection and radiation on the top surface $\partial V_{\mathrm{top}}$ is negligible compared with conduction inside the solid, hence $\partial V_{\mathrm{top}}$ is also assumed adiabatic.  
Thus, homogeneous Neumann conditions on the lateral and top boundaries are 
\begin{equation}
  \frac{\partial T}{\partial \mathbf{n}} = 0,
  \qquad \text{on }\partial V_{\mathrm{lat}}\cup\partial V_{\mathrm{top}},
  \label{eq:neumann-bc-ori}
\end{equation}
where $\mathbf{n}$ is the outward unit normal. 
At the bottom surface, which is in contact with the build platform, a Dirichlet boundary condition is applied to prescribe the temperature: 
\begin{equation}
  T = T_c,
  \qquad \text{on }\partial V_{\mathrm{bot}}.
  \label{eq:dirichlet-bc-ori}
\end{equation}

% The boundary conditions are specified as follows:
% \begin{itemize}
%   \item Homogeneous Neumann condition (adiabatic) on the lateral and top surfaces:
%     \begin{equation}\label{eq:neumann-bc}
%       % \frac{\partial T}{\partial x_i} n_i = 0,
%       \frac{\partial T}{\partial \mathbf{n}} = 0,
%       \qquad \text{on } \partial V_{\mathrm{lat}}\cup\partial V_{\mathrm{top}},
%     \end{equation}
%   \item Prescribed temperature (Dirichlet condition) on the bottom surface:
%     \begin{equation}\label{eq:dirichlet-bc}
%       T = T_c,
%       \qquad \text{on } \partial V_{\mathrm{bot}}.
%     \end{equation}
% \end{itemize}

\subsection{Temperature Decomposition}

To facilitate numerical treatment, the temperature field $T(\mathbf{x},t)$ is decomposed into an analytical component and a numerical correction:
\begin{equation}\label{eq:temperature-decomposition}
  T(\mathbf{x},t) = \tilde T(\mathbf{x},t) + \hat T(\mathbf{x},t),
\end{equation}
where $\tilde T$ is a known analytical solution of the temperature evolution due to discretized point sources, and $\hat T$ is the complementary field to be solved numerically.  
The analytical component is expressed as a sum of $K$ point-source solutions:
\begin{equation}\label{eq:sum-of-point-temp}
  \tilde T(\mathbf{x},t) = \sum_{I=1}^{K}\tilde T^{(I)}(\mathbf{x},t),
  \qquad \tau^{(K)} \le t,\quad K \le N,
\end{equation}
with each contribution given by
\begin{equation}\label{eq:semi-temperature}
  \tilde T^{(I)}(\mathbf{x},t) =
  \frac{E^{(I)}}{\rho c_p \left[ 4 \pi \alpha (t-\tau^{(I)}) \right]^{3/2}}
  \exp\!\left(
    -\frac{\|\mathbf{x}-\mathbf{x}^{(I)}\|^2}{4 \alpha \left( t-\tau^{(I)} \right) }
  \right) H \left( t-\tau^{(I)} \right),
\end{equation}
where $\mathbf{x}^{(I)}$ is the location of the $I$-th heat source, $R^{(I)} = \| \mathbf{x}-\mathbf{x}^{(I)} \|$ is the Euclidean distance to the source, $E^{(I)}=A\, P \,\Delta t$ is the deposited energy, and $H$ is the Heaviside step function ensuring causality, the factor $A$ accounts for the laser absorptivity of metal powder. The modified time $\tau^{(I)} = t^{(I)}-r^2/8\alpha$ is introduced to avoid singularity at the point-source activation time $t^{(I)}$, and accounts for the finite radius of the laser source.

Substituting \cref{eq:temperature-decomposition} into \cref{eq:linear-heat} and isolating $\hat{T}$:
\begin{equation}\label{eq:that-heat}
  \rho c_p \frac{\partial \hat T}{\partial t}
  = k \nabla^{2}\hat T, \qquad \mathbf{x}\in V.
\end{equation}

The boundary conditions for $\hat T$ become
\begin{align}
  q_N := -k\,\frac{\partial \hat T}{\partial \mathbf{n}}
    &= k\,\frac{\partial \tilde T}{\partial \mathbf{n}},
    \quad \text{on }\partial V_{\mathrm{lat}},
    \label{eq:that-bc-n}\\
  \hat T &= T_c - \tilde T,
    \quad \text{on }\partial V_{\mathrm{bot}},
    \label{eq:that-bc-d}
\end{align}
with the initial condition $\hat T(\mathbf{x},0)=T_c$ due to $\tilde T(\mathbf{x},0)=0$.
Neumann conditions are not explicitly imposed on $\partial V_{\mathrm{top}}$ in \cref{eq:that-bc-n}, as the specified Neumann condition, i.e., the adiabatic boundary condition on the top surface, is inherently satisfied and embedded in \cref{eq:semi-temperature}, owing to the coincidence of the top surface of the build part with the origin of the semi-infinite domain.

% We now derive the weak formulation of Eq.~\eqref{eq:that-heat}. Let $w \in W_0 := \{ w \in H^1(V) \mid w = T_c - \tilde{T} \text{ on } \partial V_\mathrm{bot} \}$ be a test function. Then the weak form reads:
% \begin{equation}
%     \int_V \rho c_p \frac{\partial \hat{T}}{\partial t} w \, \mathrm{d}V = 
%     \int_V k \nabla^2 \hat{T} \cdot w \, \mathrm{d}V.
% \end{equation}

% Applying the divergence theorem to the right-hand side yields:
% \begin{equation} \label{eq:weak-form}
%     \int_V \rho c_p \frac{\partial \hat{T}}{\partial t} w \, \mathrm{d}V + 
%     \int_V k \nabla \hat{T} \cdot \nabla w \, \mathrm{d}V =
%     \int_{\partial V_\mathrm{lat}} -q_N w \, \mathrm{d}s.
% \end{equation}

% Implementation Notes:
% \begin{itemize}
%     \item Define Neumann boundary function $g_N = k \frac{\partial \Tilde{T}}{\partial x_i} n_i$, then use \verb|bcInfo.addCondition()| to impose Neumann conditions on the four lateral boundaries.
%     \item {\color{red}Confirm whether Neumann boundaries in G+Smo can be time-dependent.}
%     \item Assemble the right-hand-side vector $f$ incorporating the boundary integral.
% \end{itemize}

% \begin{itemize}
%     \item Heat equation
%     \item Decomposition approach
%     \item Boundary conditions
% \end{itemize}

\begin{figure}[H]
  \centering
  \includegraphics[width=0.9\linewidth]{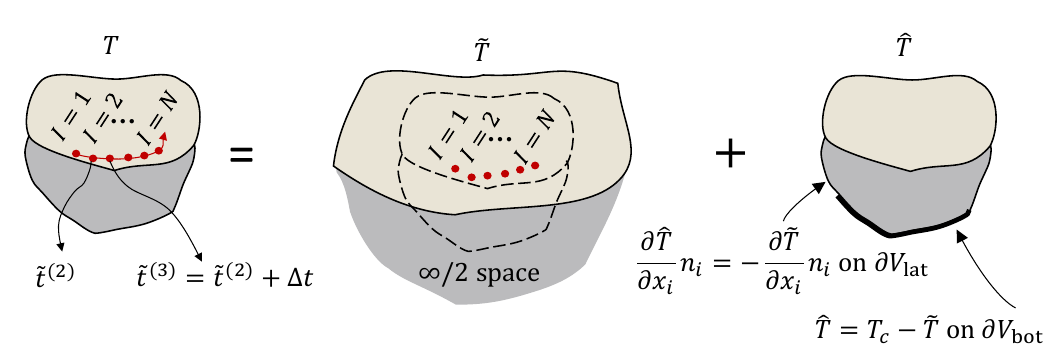}
  \caption{A curvilinear (contour) laser scan discretized into a sequence of point sources. The total temperature field is obtained by superimposing the semi-infinite point-source solutions with the complementary field that enforces the finite-part boundary conditions.}
  \label{fig:curvilinear-scan}
\end{figure}

The temperature field is decomposed as $T = \Tilde{T} + \hat{T}$ according to \eqref{eq:temperature-decomposition}. 
The analytical component $\Tilde{T}$ represents the superposition of point-source contributions of the discretized laser scan (\cref{fig:curvilinear-scan}), while the correction field $\hat{T}$ enforces the adiabatic and Dirichlet boundary conditions of the solid domain. 
We solve numerically for $\hat{T}$ using an isogeometric Galerkin method combined with an implicit time-integration scheme.

\subsection{Weak form}
\label{subsec:weak-form}

% Substituting \eqref{eq:temperature-decomposition} into \eqref{eq:heat-conduction} and isolating $\hat{T}$ yields

The strong form of the governing equation for the $\hat{T}$ field is expressed in \cref{eq:that-heat}, with the modified Neumann and Dirichlet boundary conditions specified in \cref{eq:that-bc-n,eq:that-bc-d}.
Let $w \in H^1(V)$ be a test function vanishing on $\partial V_\mathrm{bot}$ after Dirichlet enforcement. 
Multiplying the governing equation by $w$ and integrating over $V$, and applying the divergence theorem, we obtain
\begin{equation}\label{eq:weak-form-hat}
    \int_V \rho c_p \frac{\partial \hat{T}}{\partial t} w \, dV 
    + \int_V k \nabla \hat{T} \cdot \nabla w \, dV
    = \int_{\partial V_\mathrm{lat}} (-q_N) w \, dS ,
\end{equation}
with $q_N$ defined in \eqref{eq:that-bc-n}.  
This weak form incorporates the laser-induced heat flux implicitly through $q_N$ and preserves the adiabatic boundaries on $\partial V_\mathrm{lat}$.

\subsection{Isogeometric discretization}
\label{subsec:isogeometric discretization}

The geometry of the printed part is exactly represented by a NURBS mapping $\boldsymbol{\Phi}:\widehat{V}\!\to\!V$ from the parametric domain $\widehat{V}$.  
Let $\{N_i(\boldsymbol{\xi})\}_{i=1}^n$ be the B-spline or NURBS basis functions of order $p$ defined on the knot vectors of $\widehat{V}$.  
The unknown field $\hat{T}$ and the test function $w$ are approximated as
\begin{equation}
    \hat{T}_h(\boldsymbol{\xi},t) = \sum_{i=1}^{n} N_i(\boldsymbol{\xi}) \hat{T}_i(t), 
    \qquad 
    w_h(\boldsymbol{\xi}) = N_j(\boldsymbol{\xi}) .
\end{equation}
Using the standard isoparametric concept,
\begin{equation}
    \nabla_{\mathbf{x}} N_i = \mathbf{J}^{-\top} \nabla_{\boldsymbol{\xi}} N_i, \qquad dV = J \, d\boldsymbol{\xi},
\end{equation}
with $\mathbf{J} = \partial \mathbf{x}/\partial \boldsymbol{\xi}$, $\nabla_{\mathbf{x}}= \left( \partial/\partial x_1, \partial/\partial x_2, \partial/\partial x_3 \right)^{\top}$, $\nabla_{\boldsymbol{\xi}}= \left( \partial/\partial \boldsymbol{\xi}_1, \partial/\partial \boldsymbol{\xi}_2, \partial/\partial \boldsymbol{\xi}_3 \right)^{\top}$, and $J=\det(\mathbf{J})$. 
Substituting into \eqref{eq:weak-form-hat} and assembling element contributions lead to the semi-discrete matrix system
\begin{equation}
    \mathbf{M}\, \dot{\mathbf{\hat{T}}}(t) + \mathbf{K}\, \mathbf{\hat{T}}(t) = \mathbf{F}(t),
\end{equation}
where
\begin{align}
    M_{ij} &= \int_{\widehat V} \rho c_p N_i N_j J \, d\boldsymbol{\xi}, \\
    K_{ij} &= \int_{\widehat V} k (\nabla_{\mathbf{x}} N_i) \cdot (\nabla_{\mathbf{x}} N_j) J \, d\boldsymbol{\xi}, \\
    F_i(t) &= \int_{\partial V_\mathrm{lat}} (-q_N) N_i \, dS .
\end{align}
Dirichlet conditions \eqref{eq:that-bc-d} are enforced on control variables whose basis functions have support on $\partial V_\mathrm{bot}$.  
The analytical component $\Tilde{T}$ is updated explicitly at each laser position and contributes only to the boundary terms via $q_N$.

\subsection{Time integration}
\label{subsec:time integration}

The semi-discrete system is integrated using an implicit $\theta$-method:
\begin{equation}\label{eq:theta-method}
    \mathbf{M}\frac{\mathbf{\hat{T}}^{\,n+1}-\mathbf{\hat{T}}^{\,n}}{\Delta t}
    + \theta \mathbf{K}\mathbf{\hat{T}}^{\,n+1}
    + (1-\theta)\mathbf{K}\mathbf{\hat{T}}^{\,n}
    = \theta \mathbf{F}^{\,n+1} + (1-\theta)\mathbf{F}^{\,n},
\end{equation}
where $0\leq \theta \leq 1$, and the superscripts $n$ and $n+1$ refer to the current and next time steps, respectively. 
We adopt $\theta=0.5$ (Crank–Nicolson) for second-order accuracy and unconditional stability, consistent with the temporal resolution $\Delta t=\SI{1e-5}{s}$ used in the laser discretization (\cref{sc: Mathematical Model}).  
At each time step, the linear system
\begin{equation}
    \left(\frac{\mathbf{M}}{\Delta t} + \theta \mathbf{K}\right)\mathbf{\hat{T}}^{\,n+1}
    = \left(\frac{\mathbf{M}}{\Delta t} - (1-\theta)\mathbf{K}\right)\mathbf{\hat{T}}^{\,n}
    + \theta \mathbf{F}^{\,n+1} + (1-\theta)\mathbf{F}^{\,n}
\end{equation}
is solved for computing $\mathbf{\hat{T}}^{\,n+1}$.  
The reconstructed temperature field is then obtained as
\begin{equation}
    T^{\,n+1} = \Tilde{T}^{\,n+1} + \hat{T}^{\,n+1},
\end{equation}
which accounts for the transient laser heating and satisfies the boundary conditions of the printed part.

\section{Numerical examples} 
\label{sc: numerical examples}

To evaluate the accuracy and computational efficiency of the proposed semi-analytical IGA method, we present a sequence of numerical experiments. The cases progress from a single point heat source (\cref{sc:single-heat-source}), to continuous laser scanning along a curved boundary (\cref{sc:continuous-scan-boundary}), and finally culminate in a complex freeform geometry (\cref{sc:continuous-scan-complex-geometry}).

\subsection{Single heat source comparison} 
\label{sc:single-heat-source}

\begin{figure}[H]
  \centering
  \includegraphics[width=0.8\linewidth]{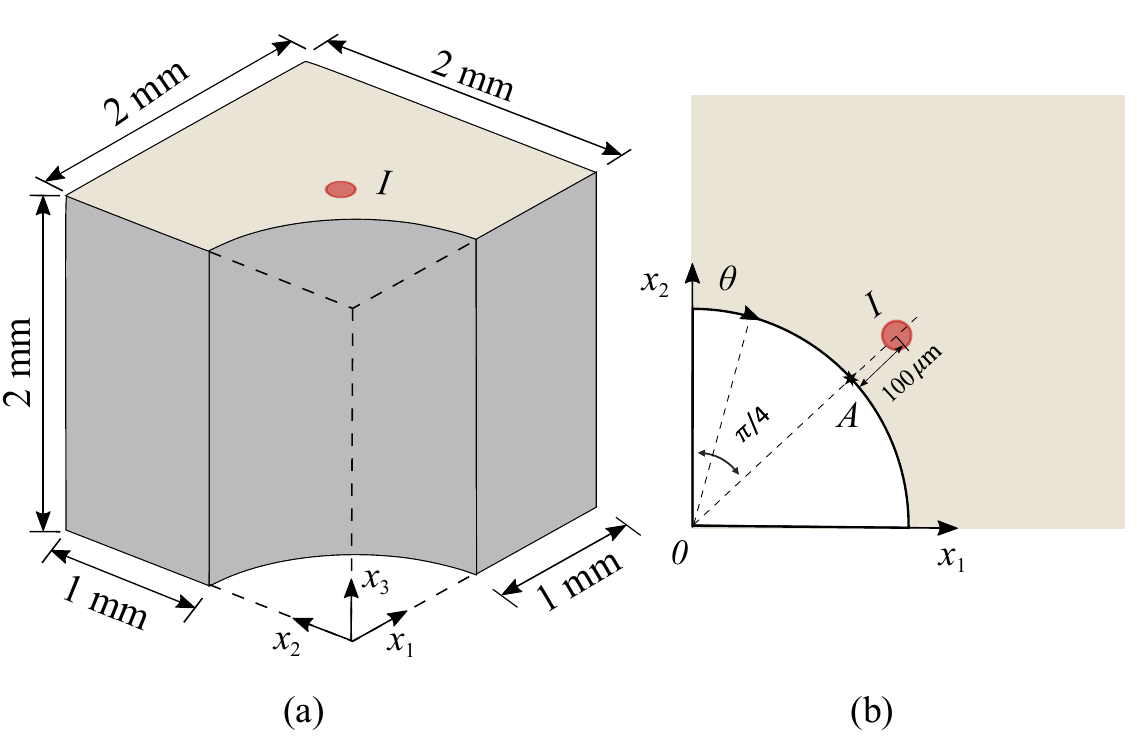}
  \caption{(a) A point source (red) positioned near the curved boundary of a 
  $\SI{2}{mm}\times\SI{2}{mm}\times\SI{2}{mm}$ cubic domain, where a quarter-cylinder 
  of radius $R_c=\SI{1}{mm}$ has been removed.  
  (b) Top view of the part. The shortest distance from the point source to the arc 
  is $\SI{100}{\micro\meter}$, located at $\theta=\pi/4$ relative to the $x_2$-axis.}
  \label{fig:single-source-position}
\end{figure}

In this section, we compare the simulation efficiency and predictive accuracy of FEM and IGA using a single point heat source. The test geometry is an extruded part obtained by removing a quarter-cylinder from a cube (\cref{fig:single-source-position}); the cube measures $\SI{2}{mm}\times\SI{2}{mm}\times\SI{2}{mm}$ and the cylindrical radius is $R_c=\SI{1}{mm}$.

A laser pulse of duration $\Delta t=\SI{1.0e-5}{s}$ and power $P=\SI{82.5}{W}$ is applied to the top surface. Because the laser spot radius is only $r_\mathrm{laser} = \SI{20}{\micro\meter}$ and the exposure time is extremely short, the heat input is approximated as a point source $I$ (red circle in \cref{fig:single-source-position}). The point source is positioned $\SI{100}{\micro\meter}$ from the curved boundary, with the line connecting it to the coordinate origin forming a $\pi/4$ angle with the $x_2$-axis (\cref{fig:single-source-position}b). The metal powder absorptivity is $A=0.77$, giving an applied thermal load of each point source $E = A\,P\Delta t= \SI{6.35E-4}{\joule}$. The thermal properties of Ti–6Al–4V used in all simulations are listed in \cref{tb:ti6al4v-thermal-properties}. 

\begin{table}[H]
  \caption{Thermal properties of Ti–6Al–4V.}
  \label{tb:ti6al4v-thermal-properties}
  \centering
  \begin{tabular}{ccc}
    \toprule
    Conductivity & Heat capacity & Density \\
    \midrule
    \SI{42}{W/(m.K)} & \SI{990}{J/(kg.K)} & \SI{4420}{kg/m^3} \\
    \bottomrule
  \end{tabular}
\end{table}

Both FEM and IGA are employed to compute the complementary temperature field $\hat{T}$ defined in \cref{eq:that-heat,eq:that-bc-n,eq:that-bc-d}, induced by the point source $I$ shown in \cref{fig:single-source-position}.

The control points of the part used for the NURBS representation are given in \cref{tb: quarter-cylinder part}, and knot vectors in three directions are shown as 
$\Xi = \{0, 0, 0, 0.5, 0.5, 1, 1, 1\}, \quad H = \{0, 0, 1, 1\}, \quad Z = \{0, 0, 1, 1\}.$

\begin{table}[H]
  \caption{Control net for the top surface of the part (unit: \SI{}{mm}).}
  \label{tb: quarter-cylinder part}
  \centering
  \begin{tabular}{ccc}
    \toprule
    $i$ & $B_{i,1}$ & $B_{i,2}$ \\
    \midrule
    1 & $(1, 0, 2)$ & $(0, 0, 2)$ \\
    2 & $(1, \sqrt{2}-1, 2)$ & $(0, 1, 2)$ \\
    3 & $(2-\sqrt{2}/2, \sqrt{2}/2, 2)$ & $(0, 2, 2)$ \\
    4 & $(3-\sqrt{2}, 1, 2)$ & $(1, 2, 2)$ \\
    5 & $(2, 1, 2)$ &  $(2, 2, 2)$ \\
    \bottomrule
  \end{tabular}
\end{table}

\begin{figure}[H]
  \centering
  \subfigure[IGA mesh, minimum element size around \SI{50}{\micro\meter}, 4845 DOFs]{\includegraphics[width=0.45\linewidth]{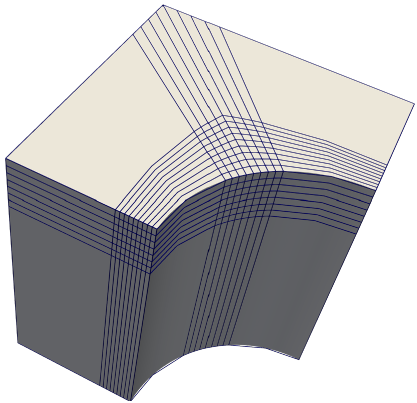}} \ 
  \subfigure[FE mesh, minimum element size around \SI{50}{\micro\meter}, 3302 DOFs]{\includegraphics[width=0.45\linewidth]{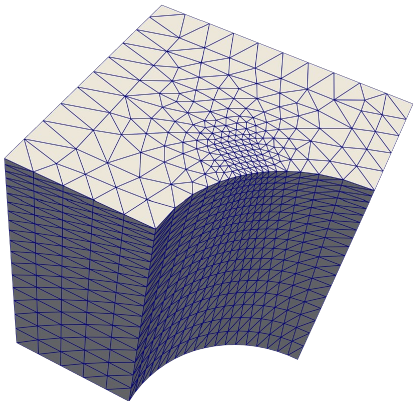}}
  \caption{IGA mesh and Finite Element (FE) mesh for computing $\hat{T}$. Small elements of approximately \SI{50}{\micro\meter} are used to refine the boundary near the point source.}
  \label{fig:mesh-iga-fem}
\end{figure}

For IGA, meshes of varying mesh density are employed to calculate the $\hat{T}$ field. One representative IGA mesh is shown in \cref{fig:mesh-iga-fem}a, where the minimum element size is approximately \SI{50}{\micro\meter}.

For the Finite Element (FE) simulations, a series of meshes with varying mesh densities, from coarse to highly refined, is generated and analyzed using ABAQUS. An example of a FE mesh is provided in \cref{fig:mesh-iga-fem}b. In this mesh, size of elements near the point source are refined to approximately \SI{50}{\micro\meter} to resolve steep temperature gradients accurately, while coarser elements of about \SI{200}{\micro\meter} are used farther apart from the source position to reduce computational cost. Additionally, a highly refined mesh with a minimum element size of \SI{5}{\micro\meter} near the source is adopted and designated as the reference solution for accuracy assessment.  

To evaluate and compare the simulation efficiency and predictive accuracy of FEM and IGA, we examine the temperature at point $A$ (marked in \cref{fig:single-source-position}b). The source is activated at $t=\SI{0}{}$, and the comparison is performed at $t=\SI{1.9e-4}{s}$, which corresponds to the time instance when the complementary temperature $\hat{T}$ at point $A$ reaches its peak value.
At this moment, the temperature gradient at the point is nearly at its maximum magnitude.
The reference FE simulation with a minimum element size of \SI{5}{\micro\meter} yields $\hat{T}_\mathrm{ref}=\SI{20.79}{^\circ\mathrm{C}}$, which is adopted to benchmark accuracy.

\begin{figure}[H]
    \hfill
    \subfigure[]{
    \includegraphics[width=0.46\linewidth]{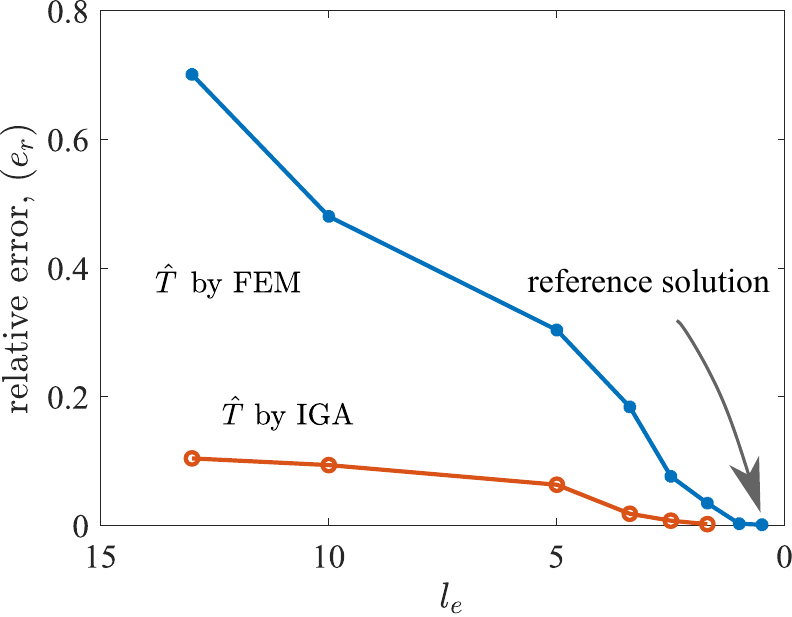}
    \label{fg:single_source_error_mesh}
    }
    \hfill
    \subfigure[]{
    \includegraphics[width=0.46\linewidth]{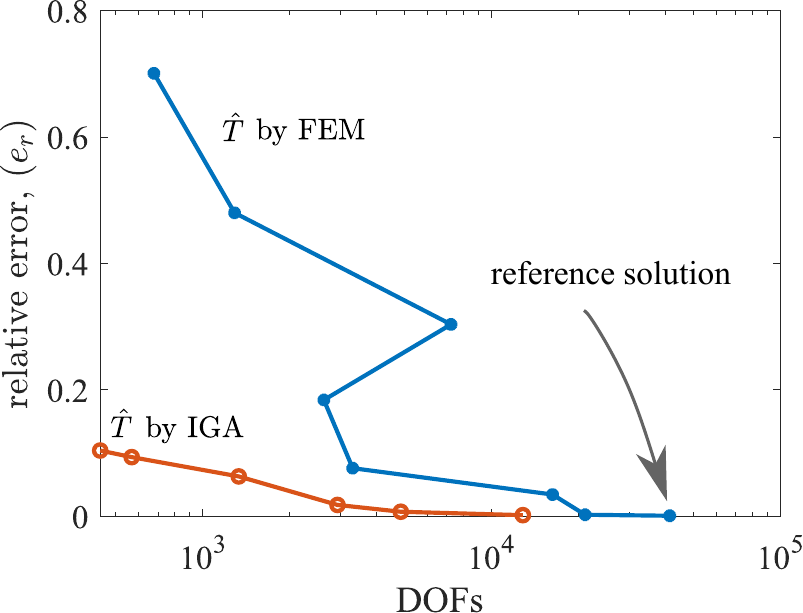}
    \label{fg:single_source_error_freedom}
    }
    \caption{Relative error $e_r$ of $\hat{T}$ solved by FEM and IGA with respect to (a) $l_e$, the ratio of minimum element size ($l_\mathrm{min}$) to laser spot radius ($r_{laser}$) and (b) number of numerical degrees of freedom. The reference result is calculated using FEM with linear tetrahedron elements of minimum element size around $\SI{5}{\micro\meter}$ (C3D4 element in ABAQUS).} \label{fg:single_source_error}
\end{figure}

As the mesh is coarsened, the number of degrees of freedom (DOFs) decreases, and the computed $\hat{T}$ field at point $A$ gradually deviates from the reference, introducing a relative error, calculated as 
\begin{equation*}
    e_r= \frac{ |\hat{T} - \hat{T}_\mathrm{ref}|}{ \hat{T}_\mathrm{ref} }.
\end{equation*}
\Cref{fg:single_source_error}a reports the relative error ($e_r$) as a function of the ratio ($l_e$) between the minimum element size ($l_\mathrm{min}$) and the laser spot radius ($r_\mathrm{laser}$)
\begin{equation*}
    l_e= \frac{l_\mathrm{min}}{r_\mathrm{laser}}.
\end{equation*}
\cref{fg:single_source_error}b illustrates the variation of $e_r$ with respect to DOFs.

In this study, the minimum element size $l_\mathrm{min}$ is normalized by the laser spot radius. The laser spot radius serves as a key parameter that characterizes the spatial intensity of the energy distribution, as shown by $\tau^{(I)}$ in \cref{eq:semi-temperature}. Since accurately capturing the steep temperature gradients induced by the laser requires the mesh resolution to be on the order of the laser spot radius \cite{ROBERTS2009916,FOROOZMEHR2016255,VANINI2024104369,ZHANG2023108839}, this normalization provides a consistent measure of discretization quality, as discussed in \cref{sec:introduction}. Although the heat flux intensity decreases as energy conducts away from the source, leading to smaller thermal gradients if the source-to-boundary distance increases. 
Therefore, the minimum element size could, in principle, also be normalized by this distance. For the boundary value problem to solve the complementary field $\hat{T}$, this is chosen by us as $\SI{100}{\micro\meter}$ in this paper. We know that heat sources distant from the boundary are not a problem because they smooth out. However, we will inevitably have sources close to the boundary since the entire top surface of the part should be scanned. And, the nearest distance from the heat source to the boundary is influenced by both the laser spot radius and the input laser energy. To maintain a consistent and general analysis, the minimum element size is therefore normalized by the laser spot radius.

\begin{figure}[H]
    \centering
    \includegraphics[width=1\linewidth]{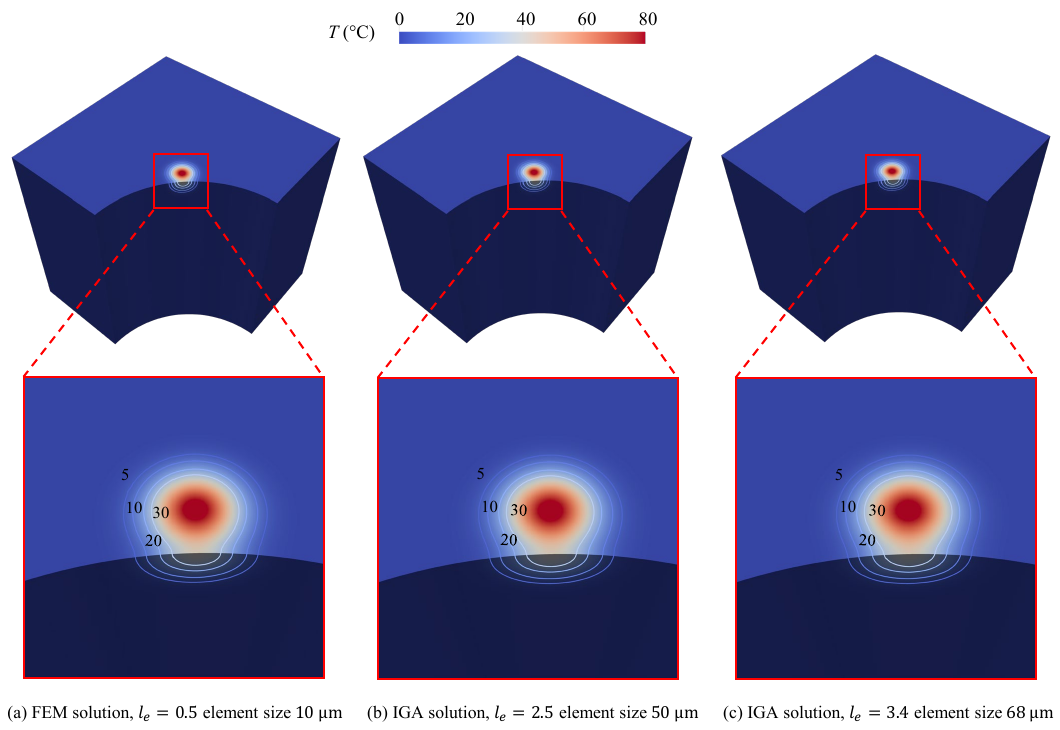}
    \caption{Temperature distribution at $t = \SI{1.9E-4}{s}$. Three simulations with different numerical methods and minimum element sizes are adopted. }
    \label{fig:temperature_single_source}
\end{figure}

% According to \cref{fg:single_source_error}, the relative error $e_r$ increases drastically when the minimum element size $l_e = 2.5$ for FEM, implying a minimum element size of \SI{50}{\micro\meter}, 2.5 times of the laser spot radius, and corresponding to 3302 DOFs. Interestingly, when $l_e$ increases from 3.4 to 5, the DOF count unexpectedly rises from 2624 to 7231. This counterintuitive behavior results from the meshing strategy: for $l_e=3.4$, i.e., the minimum element size \SI{68}{\micro\meter}, local refinement is applied near the laser spot while coarser elements are used elsewhere, producing a relatively low DOF count. In contrast, the minimum element size \SI{100}{\micro\meter} ($l_e=5$) is too large for the small overall dimensions (\SI{2}{\milli\meter} cube), thus removing the need for local refinement. Therefore, a uniform mesh is adopted, substantially increasing the total DOFs despite the larger element size. 

% For the IGA simulations, the relative error $e_r$ is evaluated against the FEM reference solution with a minimum element size of $\SI{5}{\micro\meter}$ ($l_e=0.25$). As shown in \cref{fg:single_source_error}a, $e_r$ remains around 0.1 (10\%) even when the IGA element size is as large as $l_e=13$. This comparison indicates that, for a given error tolerance, IGA requires fewer degrees of freedom and permits substantially larger element sizes than FEM.

According to \cref{fg:single_source_error}, the relative error $e_r$ for FEM increases drastically when the non-dimensional minimum element size reaches $l_e=2.5$, which corresponds to a physical minimum element size of \SI{50}{\micro\meter}, i.e., 2.5 times the laser spot radius, and 3302 DOFs. At first glance, it may appear counterintuitive that when $l_e$ increases from 3.4 to 5, the DOF count rises from 2624 to 7231. This behavior is a direct consequence of the meshing rules adopted in this study rather than an arbitrary tuning of the discretization.

For FEM, we employ the same meshing strategy for all values of $l_e$: the target minimum element size is prescribed by $l_e$, and (i) local refinement is applied in a small region around the laser spot, while (ii) away from the source the mesh is kept as coarse as permitted by a mesh-quality criterion (a minimum number of elements across each spatial direction). For $l_e=3.4$ (minimum element size \SI{68}{\micro\meter}), this results in a strongly graded mesh with local refinement near the heat source and coarse elements elsewhere, yielding only 2624 DOFs. When $l_e$ is further increased to 5 (minimum element size $\SI{100}{\micro\meter}$) in a cube of edge length $\SI{2}{\milli\meter}$, this grading strategy would produce fewer than the required number of elements in the coarser regions and thus violate the mesh-quality constraint. In this regime the procedure naturally degenerates into an (approximately) uniform mesh, which explains why the DOF count increases despite the larger element size. We stress that the mesh is not manually adjusted per data point; instead, the DOFs follow from a fixed, physics- and quality-driven meshing procedure as $l_e$ varies.

For the IGA simulations, the relative error $e_r$ is evaluated against the FEM reference solution with a minimum element size of $\SI{5} {\micro\meter}$ ($l_e=0.25$). As shown in \cref{fg:single_source_error}a, $e_r$ remains around 0.1 (10\%) even when the IGA element size is as large as $l_e = 13$. Despite the non-monotonic relationship between $l_e$ and the DOF count on the FEM side, this comparison is still objective: each data point for FEM and IGA corresponds to a mesh generated by a predefined, method-consistent meshing strategy under the same physical problem and error metric. The observed difference in DOF is thus a consequence of the underlying approximation spaces (continuous Lagrange elements versus highly continuous spline bases) rather than of any ad hoc adjustment of the meshes. Consequently, the results in \cref{fg:single_source_error} robustly indicate that, for a given error tolerance, IGA can achieve the same accuracy with fewer DOFs and substantially larger element sizes than FEM.

\begin{figure}[H]
    \hfill
    \subfigure[]{\includegraphics[width = 0.31\textwidth]{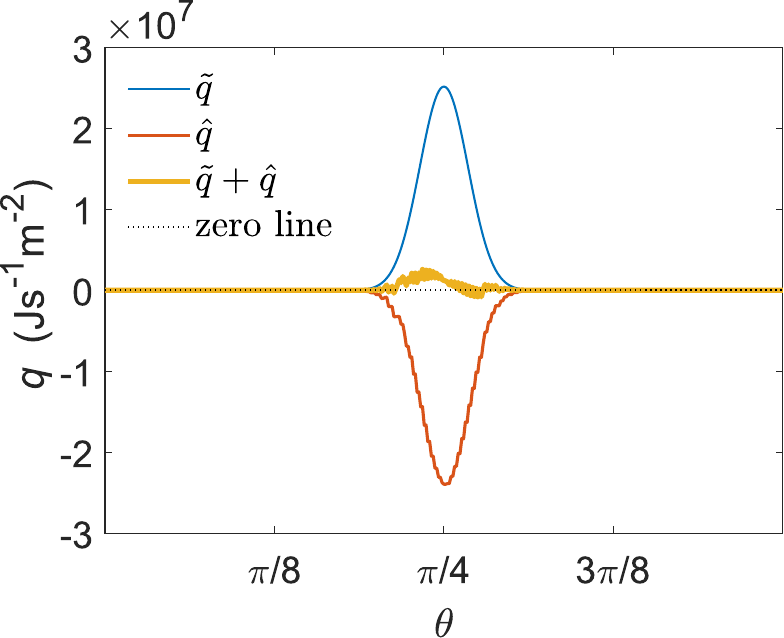}
    \label{fg:one_source_flux_FEM_10um}}
    \hfill
    \subfigure[]{\includegraphics[width = 0.31\textwidth]{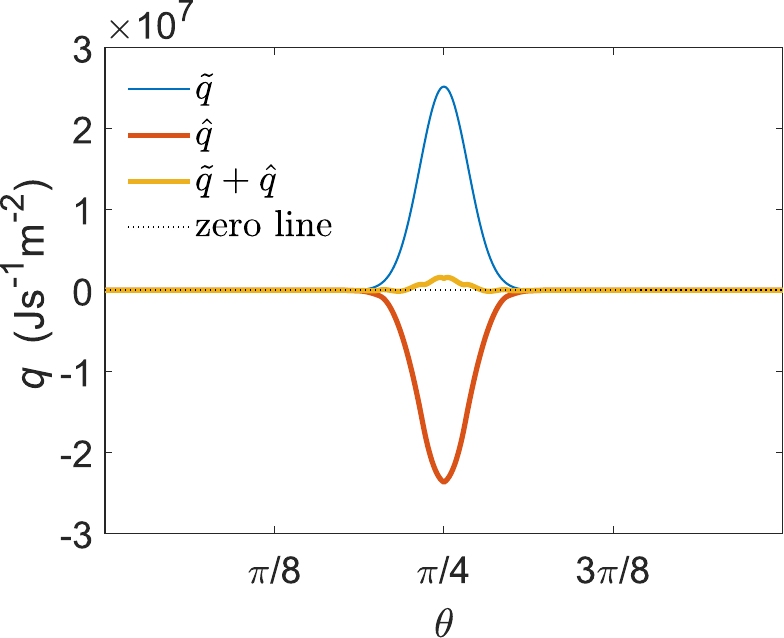}
    \label{fg:one_source_flux_IGA_50um}}  
    \hfill
    \subfigure[]{\includegraphics[width = 0.31\textwidth]{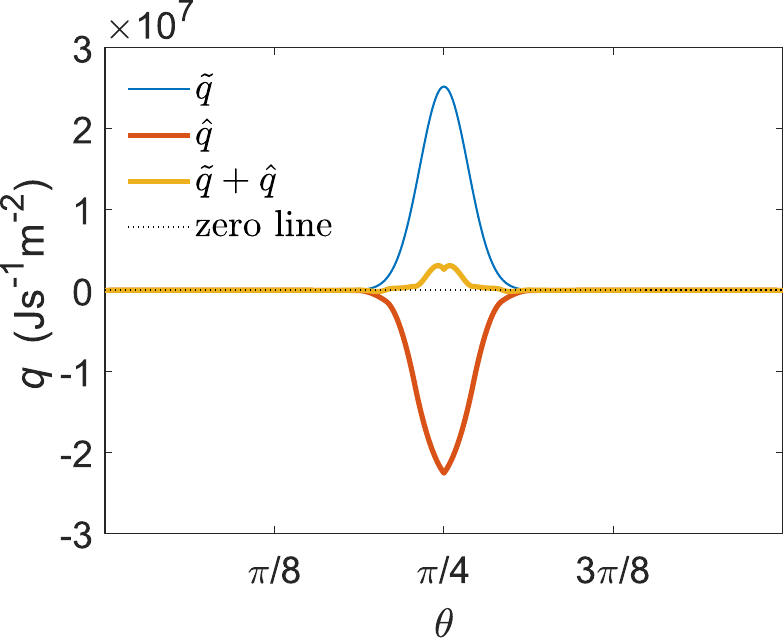}
    \label{fg:one_source_flux_IGA_68um}}
     \hfill\strut
    \caption{Heat loss rate distribution $q$ along the curve boundary on the top surface for the three simulations with different methods and minimum element sizes at $t = \SI{1.9E-4}{s}$. The angular coordinate $\theta$ is defined with respect to the $x_2$-axis shown in \cref{fig:single-source-position}. (a) $\hat{q}$ solved by FEM, $l_e = 0.50$ the minimum element size on the boundary is around \SI{10}{\micro\meter}, and $\hat{T}$ is solved using linear tetrahedron elements (C3D4 element in ABAQUS). (b) $\hat{q}$ solved by IGA, $l_e = 2.5$ the minimum element size on the boundary is around \SI{50}{\micro\meter}. (c) $\hat{q}$ solved by IGA, $l_e = 3.4$ the minimum element size on the boundary is around \SI{68}{\micro\meter}.} \label{fg:one_source_flux_comparison}
\end{figure}

The total temperature ($\tilde{T}+\hat{T}$) distribution produced by the point source at $t=\SI{1.9e-4}{s}$ is presented in \cref{fig:temperature_single_source}. Figs. \ref{fig:temperature_single_source}b and c display temperature contours obtained with the IGA method with $l_e = 2.5$ and $3.4$, respectively, while \cref{fig:temperature_single_source}a shows the corresponding FEM result with $l_e = 0.50$. Please note that we do not use the reference case $l_e = 0.25$ because the results in \cref{fg:single_source_error}b indicate that the temperature has already converged.
In all cases, the temperature isolines are orthogonal to the part boundary surfaces, confirming that the adiabatic boundary conditions are satisfied. 
Besides, we also calculate the net heat loss rate,
\begin{equation}
    \tilde{q}+\hat{q} := -k\,\frac{\partial T}{\partial \mathbf{n}} = -k\,(\frac{\partial\tilde{T}}{\partial \mathbf{n}}+\frac{\partial\hat{T}}{\partial \mathbf{n}}),
\end{equation}
along the curve boundary of the part to compare the differences between IGA and FEM. For the exact adiabatic boundary condition, the net heat loss rate should be zero everywhere along the boundary, as illustrated in \cref{eq:neumann-bc-ori}.
As shown in \cref{fg:one_source_flux_comparison}, the blue line in the three subfigures represents the heat loss rate ($\tilde{q}$) due to the point source analytical solution ($\tilde{T}$) and is identical across the three subfigures, while the orange line represents the heat loss rate ($\hat{q}$) due to the complementary field ($\hat{T}$). The yellow line is the net heat loss distribution due to the superposition of two fields ($\tilde{q}+\hat{q}$).
It can be observed that the $\hat{q}$ field compensates the $\tilde{q}$ to achieve adiabatic boundary conditions, while there is still some residual heat loss through the boundary.
The heat loss rate distribution obtained using IGA with a mesh size of $l_e = 2.5$ is lower than that of by using FEM with a mesh size of $l_e = 0.50$. Furthermore, the IGA results with a mesh size of $l_e = 3.4$ demonstrate a performance comparable to that of FEM with a mesh size of $l_e = 0.5$. 

\begin{figure}[H]
  \centering
  \includegraphics[width=0.8\linewidth]{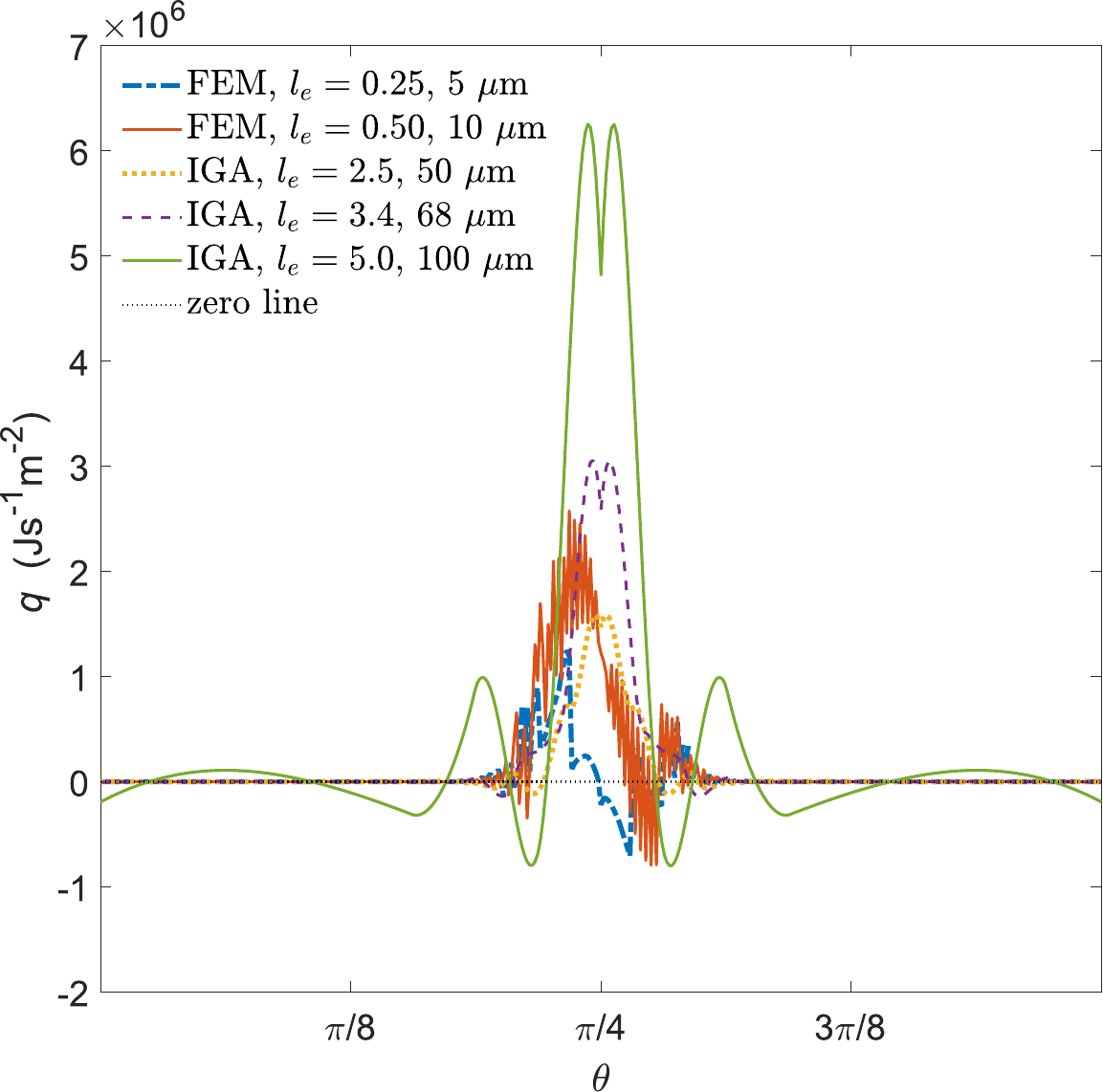}
  \caption{Net heat loss rate ($\tilde{q}+\hat{q}$) distribution along the curve boundary on the top surface at $t = \SI{1.9E-4}{s}$ from simulations with different minimum element lengths in FEM and IGA. The angular coordinate $\theta$ is defined with respect to the $x_2$-axis shown in \cref{fig:single-source-position}.}
  \label{fg:one_source_total_flux_comparison}
\end{figure}

The net heat loss rates obtained from the IGA and FE simulations with various mesh sizes are presented in \cref{fg:one_source_total_flux_comparison}. Besides, the corresponding performance metrics derived from the temperature histories are summarized in \cref{tb:metrics}, where $\int_l |\tilde{q}+\hat{q}|\mathrm{d}s$ is the integrated absolute net heat loss rate along a path $l$ and the integration path $l$ is the arc boundary on the top surface in \cref{fig:single-source-position}.
In terms of the integrated heat loss rate, the results indicate that the IGA simulations with minimum element sizes of $l_e = 2.5$ and $l_e = 3.4$ achieve high numerical accuracy, with deviations within 10\% relative to the reference case. The IGA simulations with these two mesh sizes yield results comparable to those of the FE simulation $l_e = 0.5$ (minimum element size $\SI{10}{\micro\meter}$).

\begin{table}[H]
  \caption{Metrics computed for thermal history differences at $t = \SI{1.9E-4}{s}$. The absolute difference of the maximum $\hat{T}$, $|\Delta T_\mathrm{max}|$, the integrated absolute net heat loss rate along a path $l$, $\int_l |\tilde{q}+\hat{q}|\mathrm{d}s$, the integrated absolute analytical heat loss rate along the path $l$, $\int_l |\tilde{q}|\mathrm{d}s$. The line integral path $l$ is the arc boundary on the top surface of the part in \cref{fig:single-source-position}.}
  \label{tb:metrics}
  \centering
  \begin{tabular}{cccccc}
    \toprule
      & FEM $\SI{5}{\micro\meter}$ & FEM $\SI{10}{\micro\meter}$ & IGA $\SI{50}{\micro\meter}$ & IGA $\SI{68}{\micro\meter}$ & IGA $\SI{100}{\micro\meter}$\\
      $l_e$ &  0.25 & 0.50 & 2.5 & 3.4 & 5.0 \\
    \midrule

     $|\Delta T_\mathrm{max}|=|\hat{T} - \hat{T}_\mathrm{ref}|, \, (^\circ\mathrm{C})$   & / & 0.021& 0.15 & 0.38 & 1.3 \\
      $e_r$ & / & 0.10\% & 0.73\% & 1.8\% & 6.3\% \\
    $\int_l |\tilde{q}+\hat{q}|\mathrm{d}s, \, (\mathrm{Js}^{-1}\mathrm{m}^{-1})$ & \SI{6.13E1}{} & \SI{2.62E2}{} & \SI{1.67E2}{} & \SI{3.23E2}{} & \SI{9.74E2}{} \\
     $\int_l |\tilde{q}+\hat{q}|\mathrm{d}s \Big/ \int_l |\tilde{q}|\mathrm{d}s$ & 1.74\% & 7.44\% & 4.75\% & 9.17\% & 27.64\% \\
    \bottomrule
  \end{tabular}
\end{table}

However, the IGA results with a mesh size of $l_e = 5.0$ exhibit a higher heat loss rate distribution along the boundary in  \cref{fg:one_source_total_flux_comparison,tb:metrics}. It should be noted that the single point heat source is located along the line that bisects the part into two equal halves, as illustrated in \cref{fig:single-source-position}. Along this line, according to the knot multiplicity, the solution field exhibits the lowest level of continuity, namely ($C^0$) continuity. In the subsequent example involving continuous laser scanning, the results demonstrate that the simulation using an IGA element size of $l_e = 5.0$ (minimum element size \SI{100}{\micro\meter}) achieves accuracy comparable to that obtained with an element size of $l_e = 0.5$ (minimum element size \SI{10}{\micro\meter}) in FEM.

\subsection{Continuous laser scanning along part boundary}
\label{sc:continuous-scan-boundary}

\begin{figure}[H]
  \centering
  \includegraphics[width=0.6\linewidth]{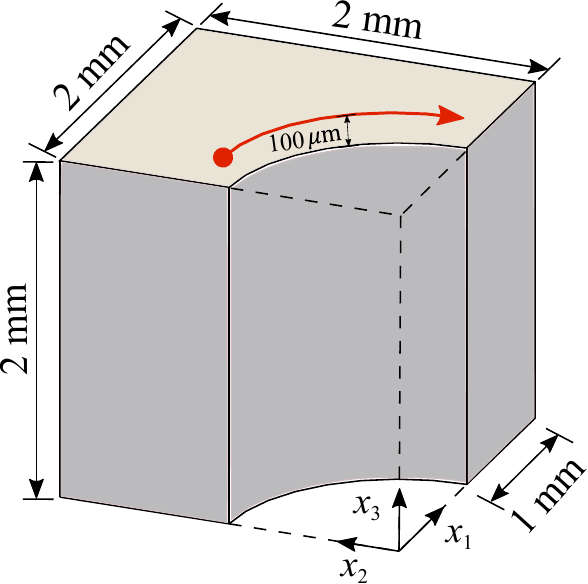}
  \caption{Geometry for the continuous contour laser scan along the part boundary. 
  The scan path is offset by \SI{100}{\micro\meter} from the boundary.}
  \label{fig:line-source-position}
\end{figure}

Building on the single point-source analysis from the previous section, we next examine a continuous laser scanning case using the same part geometry (\cref{sc:single-heat-source}, \cref{fig:single-source-position}), extruded without cross-section variation along the build ($x_3$) direction, as shown in \cref{fig:line-source-position}. 

\begin{figure}[H]
  \centering
  \subfigure[FE mesh ($l_e=0.5$, minimum element size \SI{10}{\micro\meter}), 45773 DOFs]{
  \includegraphics[width=0.48\linewidth]{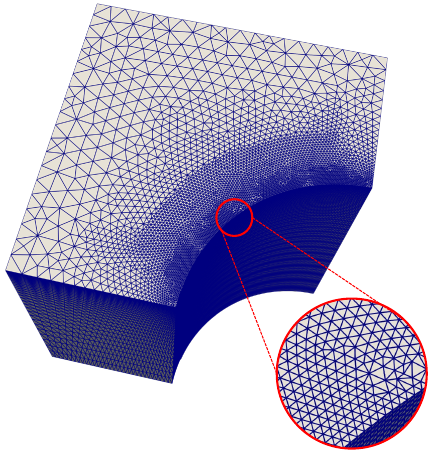}
  } \quad
  \subfigure[IGA mesh ($l_e=5.0$, minimum element size \SI{100}{\micro\meter}), 5415 DOFs]{
  \includegraphics[width=0.4\linewidth]{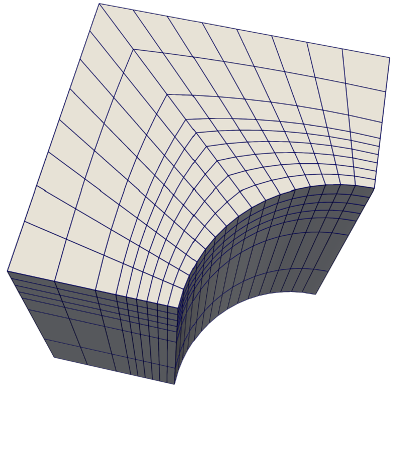}
  }
  \caption{Finite Element (FE) and isogeometric analysis (IGA) meshes 
  for the continuous laser scan simulation. 
  The left panel shows the FE mesh, and the right panel shows the IGA mesh.}
  \label{fig:line-source-mesh-comparison}
\end{figure}

The laser operates at a power of $P=\SI{82.5}{W}$, scanning speed of \SI{0.5}{\meter\per\second}, and spot radius of $r_\mathrm{laser}=\SI{20}{\micro\meter}$, with the scan path offset by \SI{100}{\micro\meter} from the part boundary. As discussed earlier, the continuous scan is discretized into a sequence of point sources using a temporal step size of $\Delta t=\SI{1e-5}{s}$. The same material properties representative of Ti-6Al-4V listed in \cref{tb:ti6al4v-thermal-properties} are used, with a metal powder absorptivity of $A=0.77$. For FEM, a minimum element size of \SI{10}{\micro\meter} ($l_e = 0.5$ finest test mesh for FEM) is adopted, resulting in 45,773 DOFs (\cref{fig:line-source-mesh-comparison}a), whereas for IGA, a minimum element size of \SI{100}{\micro\meter} ($l_e = 5.0$) with 6,137 DOFs is used (\cref{fig:line-source-mesh-comparison}b).

\begin{figure}[H]
  \centering
  \includegraphics[width=\linewidth]{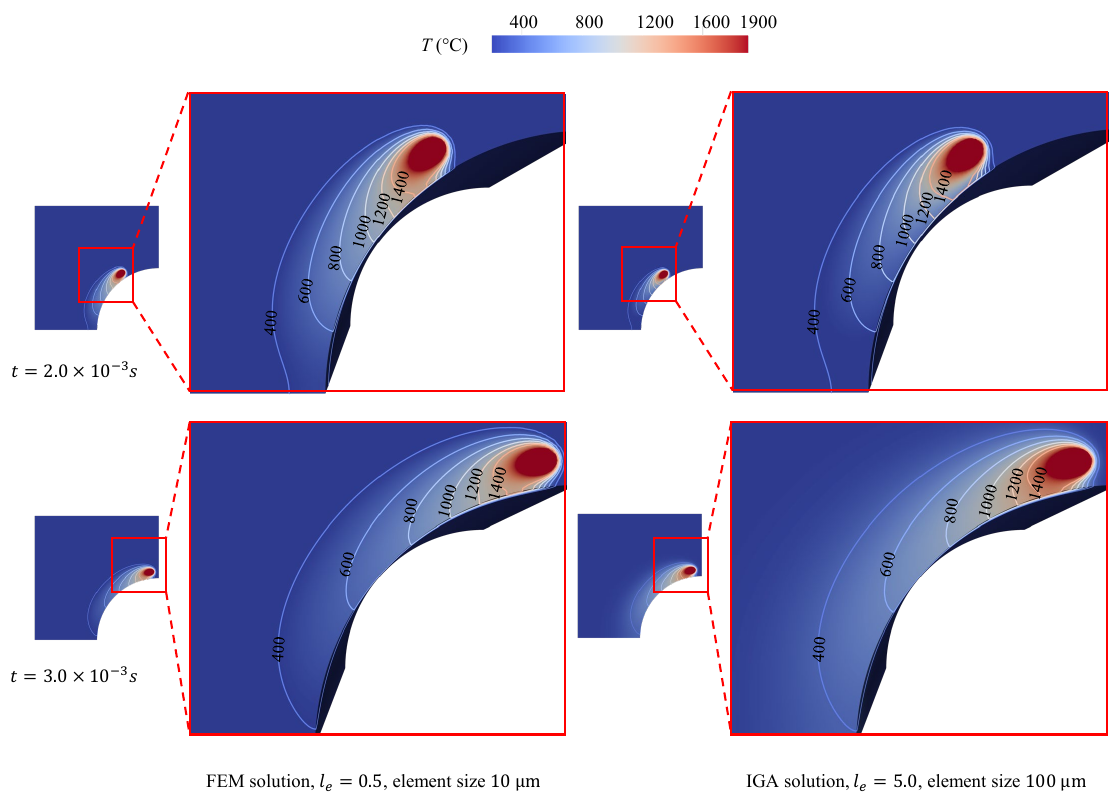}
  \caption{Temperature distributions at $t=\SI{2.0e-3}{s}$ and $t=\SI{3.0e-3}{s}$ 
  during the laser scan along the part boundary.  
  FEM uses a minimum element size of \SI{10}{\micro\meter} ($l_e = 0.5$), while IGA uses \SI{100}{\micro\meter} ($l_e = 5.0$).  
  Despite the coarser mesh, the IGA contours closely match the FEM results, 
  demonstrating comparable accuracy with significantly fewer DOFs.}
  \label{fig:single-line-source-temp}
\end{figure}

\begin{figure}[H]
    \hfill
    \subfigure[]{
    \includegraphics[width = 0.45\textwidth]{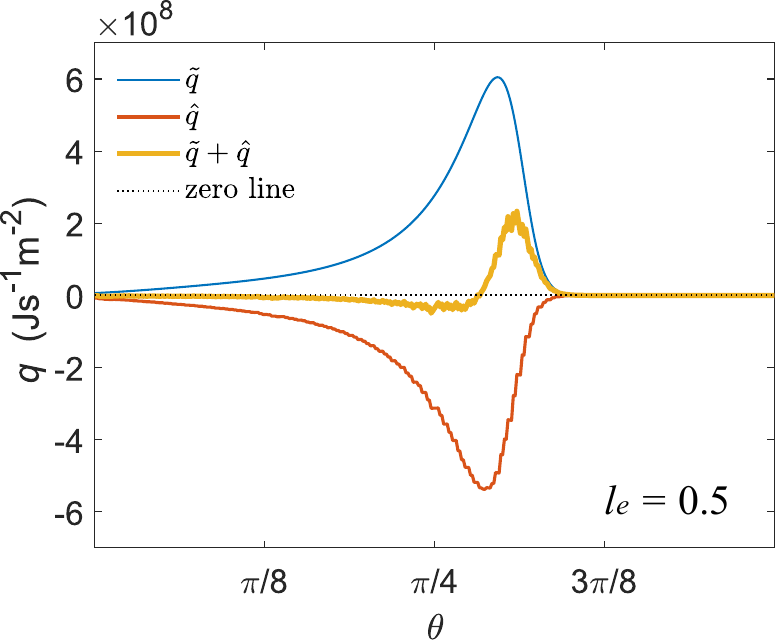}
    \label{fg:simple_line_source_FEM10um_flux_t2}
    }
    \hfill
    \subfigure[]{\includegraphics[width = 0.45\textwidth]{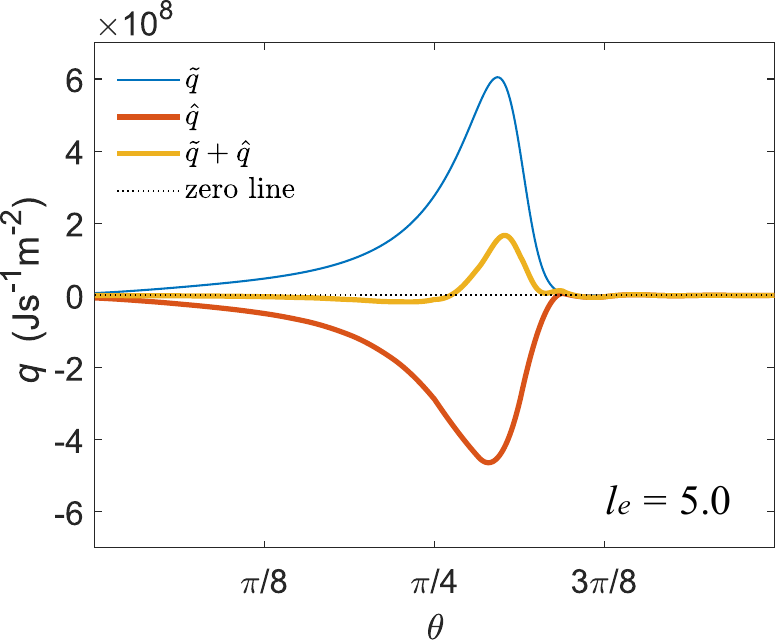}
    \label{fg:simple_line_source_IGA100um_flux_t2}}
     \hfill\strut
    \caption{Heat loss rate distribution $q$ along the curve boundary on the top surface for the two simulations with different methods and minimum element sizes at $t = \SI{2.0E-3}{s}$. The angular coordinate $\theta$ is defined with respect to the $x_2$-axis shown in \cref{fig:single-source-position}. (a) $\hat{q}$ solved by FEM with $l_e=0.5$, the minimum element size on boundary is around \SI{10}{\micro\meter}, and $\hat{T}$ is solved using linear tetrahedron elements (C3D4 element in ABAQUS). (b) $\hat{q}$ solved by IGA with $l_e=5.0$, the minimum element size on boundary is around \SI{100}{\micro\meter}.} \label{fg:simple_line_source_flux_t2}
\end{figure}

\begin{figure}[H]
    \hfill
    \subfigure[]{
    \includegraphics[width = 0.45\textwidth]{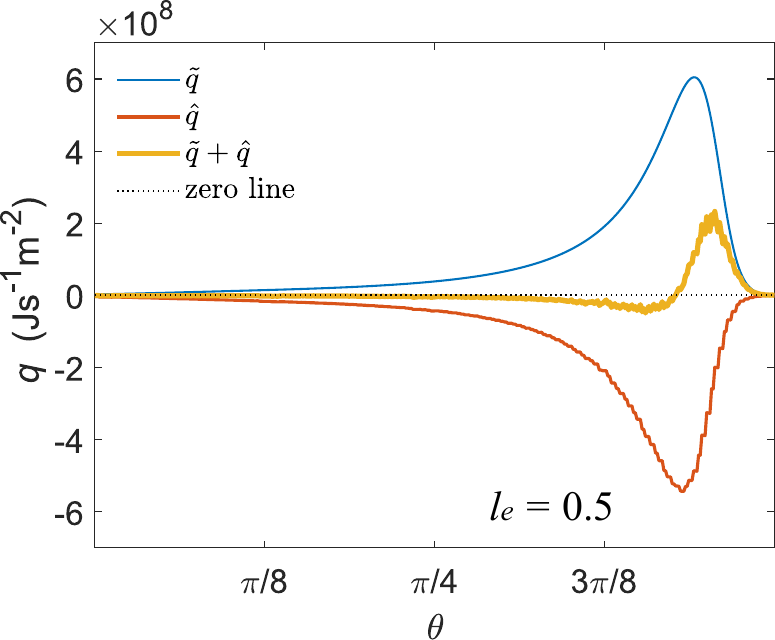}
    \label{fg:simple_line_source_FEM10um_flux_t3}
    }
    \hfill
    \subfigure[]{\includegraphics[width = 0.45\textwidth]{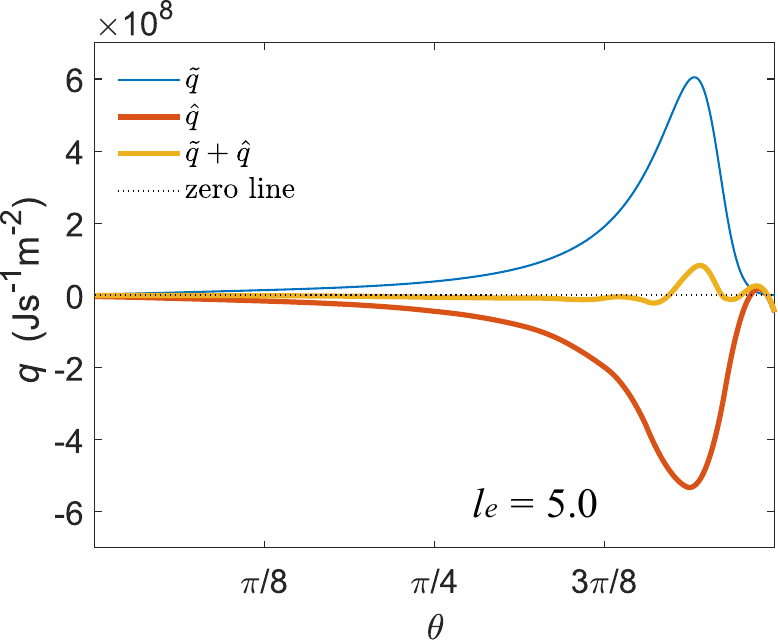}
    \label{fg:simple_line_source_IGA100um_flux_t3}}
    \hfill\strut
    \caption{Heat loss rate distribution $q$ along the curve boundary on the top surface for the two simulations with different methods and minimum element sizes at $t = \SI{3.0E-3}{s}$. The angular coordinate $\theta$ is defined with respect to the $x_2$-axis shown in \cref{fig:single-source-position}. (a) $\hat{q}$ solved by FEM $l_e=0.5$, the minimum element size on the boundary is around \SI{10}{\micro\meter}, and $\hat{T}$ is solved using linear tetrahedron elements (C3D4 element in ABAQUS). (b) $\hat{q}$ solved by IGA with $l_e=5.0$, the minimum element size on the boundary is around \SI{100}{\micro\meter}.} \label{fg:simple_line_source_flux_t3}
\end{figure}

The laser starts scanning at $t=0$, and the temperature fields at $t=\SI{2.0e-3}{s}$ and $t=\SI{3.0e-3}{s}$ obtained from FEM and IGA are shown in \cref{fig:single-line-source-temp}. 
The temperature contour lines remain orthogonal to the part boundaries, confirming that the adiabatic boundary conditions are satisfied. 
In addition, the heat loss rates along the curve boundary on the top surface obtained from the two simulation methods are compared in \cref{fg:simple_line_source_flux_t2} for the time step $t=\SI{2.0e-3}{s}$ and \cref{fg:simple_line_source_flux_t3} for the time step $t=\SI{3.0e-3}{s}$. It can be observed that at both time steps, the peak net heat loss rate on the boundary predicted by IGA with a minimum element size of $l_e = 5.0$ is lower than that obtained with the FE simulation with $l_e = 0.5$.
The corresponding integrated absolute net heat loss rates along the curve boundary are summarized in \cref{tb:metrics for continous scan}, where the line integral path $l$ is the arc boundary on the top surface of the part in \cref{fig:line-source-position}.

For the IGA simulations, the integrated absolute net heat loss rates correspond to 14.94\% and 8.03\% of the total outgoing heat loss predicted by the analytical solutions at the two respective time instances, while the corresponding values for the FE simulations are 17.76\% and 17.30\%. At $t = \SI{2.0e-3}{s}$, the numerical accuracy of IGA is slightly higher than that of FEM. In contrast, at $t = \SI{3.0e-3}{s}$, the IGA results exhibit significantly greater accuracy, as the integrated net heat loss rate is nearly half of that obtained from FEM. This improvement occurs because, at $t = \SI{2.0e-3}{s}$, the laser source is positioned closer to the $C^0$ continuity line, whereas at $t = \SI{3.0e-3}{s}$, it is located farther away, thereby reducing the effect of the lower level of continuity and enhancing the solution accuracy.

\begin{table}[H]
  \caption{Metrics computed for thermal history differences at $t = \SI{2.0E-3}{s}$ and $t = \SI{3.0E-3}{s}$.  The integrated absolute net heat loss rate along a path $l$, $\int_l |\tilde{q}+\hat{q}|\mathrm{d}s$, the integrated absolute analytical heat loss rate along the path $l$, $\int_l |\tilde{q}|\mathrm{d}s$. The line integral path $l$ is the arc boundary on the top surface of the part in \cref{fig:line-source-position}.}
  \label{tb:metrics for continous scan}
  \centering
  \begin{tabular}{c|cc|cc}
    \toprule
    & \multicolumn{2}{c|}{$t = \SI{2.0E-3}{s}$} & \multicolumn{2}{c}{$t = \SI{3.0E-3}{s}$} \\
      \midrule
    & FEM $\SI{10}{\micro\meter}$ & IGA $\SI{100}{\micro\meter}$ & FEM $\SI{10}{\micro\meter}$ & IGA $\SI{100}{\micro\meter}$\\
    $l_e$& 0.50 & 5.0 & 0.50 & 5.0\\
    \midrule
    
    $\int_l |\tilde{q}+\hat{q}|\mathrm{d}s$ & \SI{2.86E4}{} & \SI{2.41E4}{} & \SI{2.89E4}{} & \SI{1.34E4}{} \\
     $\int_l |\tilde{q}+\hat{q}|\mathrm{d}s \Big/ \int_l |\tilde{q}|\mathrm{d}s$ & 17.76\% & 14.94\% & 17.30\% & 8.03\%\\
    \bottomrule
  \end{tabular}
\end{table}

These results indicate that the overall numerical accuracy of IGA with a minimum element size of $\SI{100}{\micro\meter}$ and FEM with a minimum element size of $\SI{10}{\micro\meter}$ methods is comparable. Furthermore, the close agreement observed between the FEM and IGA temperature contour patterns demonstrates that the IGA approach, even when using a 10 times coarser mesh, can accurately reproduce the FEM results. Therefore, in the subsequent section, the same minimum element size is adopted for the IGA simulation of a more complex geometry to further evaluate its performance and computational efficiency.

\subsection{Continuous laser scanning along a complex part boundary}
\label{sc:continuous-scan-complex-geometry}

\begin{figure}[H]
  \hfill
  \subfigure[]{
  \includegraphics[width=0.47\linewidth]{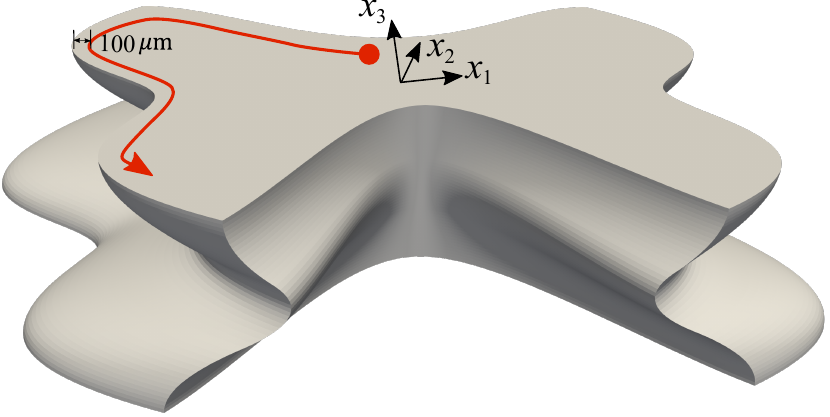}
  }
  \hfill
  \subfigure[]{
  \includegraphics[width=0.45\linewidth]{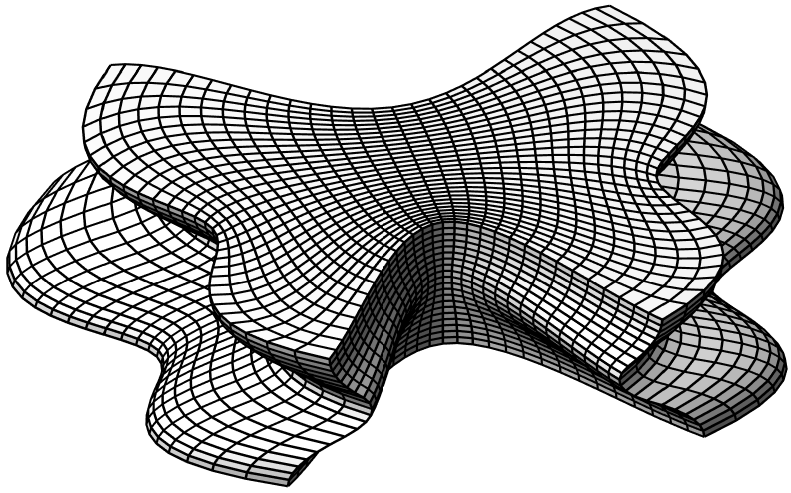}
  }
  \caption{Butterfly-shaped part constructed using NURBS and control points. The cross-sections vary nonlinearly along the build ($x_3$) direction. The complete geometric model is provided in the supplementary material. The laser contour scan is offset by \SI{100}{\micro\meter} from the boundary, and an IGA minimum element size of \SI{100}{\micro\meter} is used for the simulation.}
  \label{fig:butterfly-geo}
\end{figure}

This section extends the simulation to a more complex geometry—a butterfly-shaped part illustrated in \cref{fig:butterfly-geo}a. The geometry mapping is generated using the parameterization method of \cite{ji2021constructing}. Unlike the butterfly part in our previous work~\cite{YANG2025127059}, which maintained identical cross-sections along the build ($x_3$) direction, the present geometry features cross-sections that vary nonlinearly along the height direction. The full geometric model is provided as supplementary material.

A laser scan is performed with a \SI{100}{\micro\meter} offset from the boundary. The corresponding IGA geometry mesh is used along the scanning path shown in \cref{fig:butterfly-geo}b. A minimum element size of approximately $\SI{100}{\micro\meter}$ along the scanning path is adopted for the IGA discretization.
All process parameters—including laser power, scanning speed, spot radius, Ti–6Al–4V thermal properties, and metal powder absorptivity—are identical to those used in the previous section. This part has a length and width of approximately \SI{10}{m\meter} and a height of \qty{2}{m\meter}. The IGA simulation is conducted using the mesh with a minimum element size of \SI{100}{\micro\meter} ($l_e=5.0$), resulting in a total of \num{27456} degrees of freedom. However, to simulate the same scanning process using FEM, it leads to excessive DOFs, which is far more computationally prohibitive.

\begin{figure}[H]
  \centering
  \includegraphics[width = \linewidth]{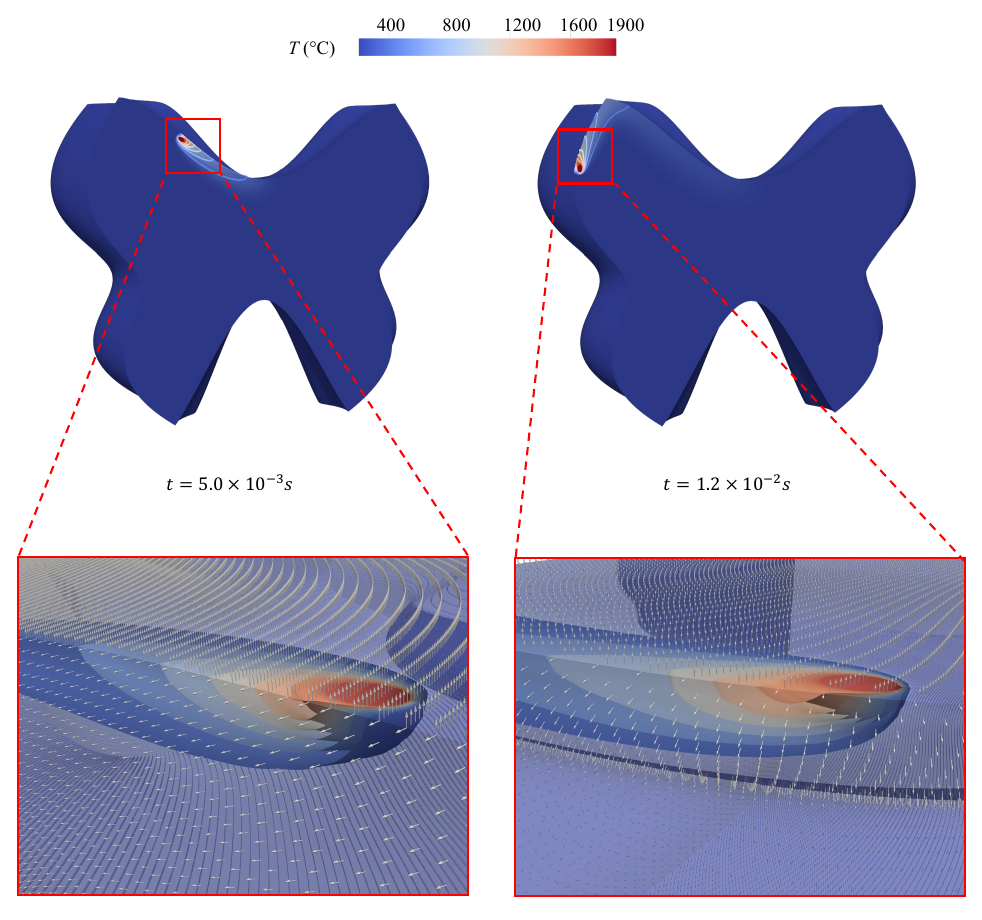}
  \caption{IGA temperature field snapshots for the butterfly-shaped part at $t=\SI{5.0e-3}{s}$ and $t=\SI{1.2e-2}{s}$ during the laser contour scan. The temperature contour surfaces remain orthogonal to the boundary surfaces, confirming proper enforcement of adiabatic boundary conditions. The white arrow represents the direction of the outward unit normal.}
  \label{fig:butterfly-temperature}
\end{figure}

\begin{figure}[H]
  \centering
  \includegraphics[width=0.88\linewidth]{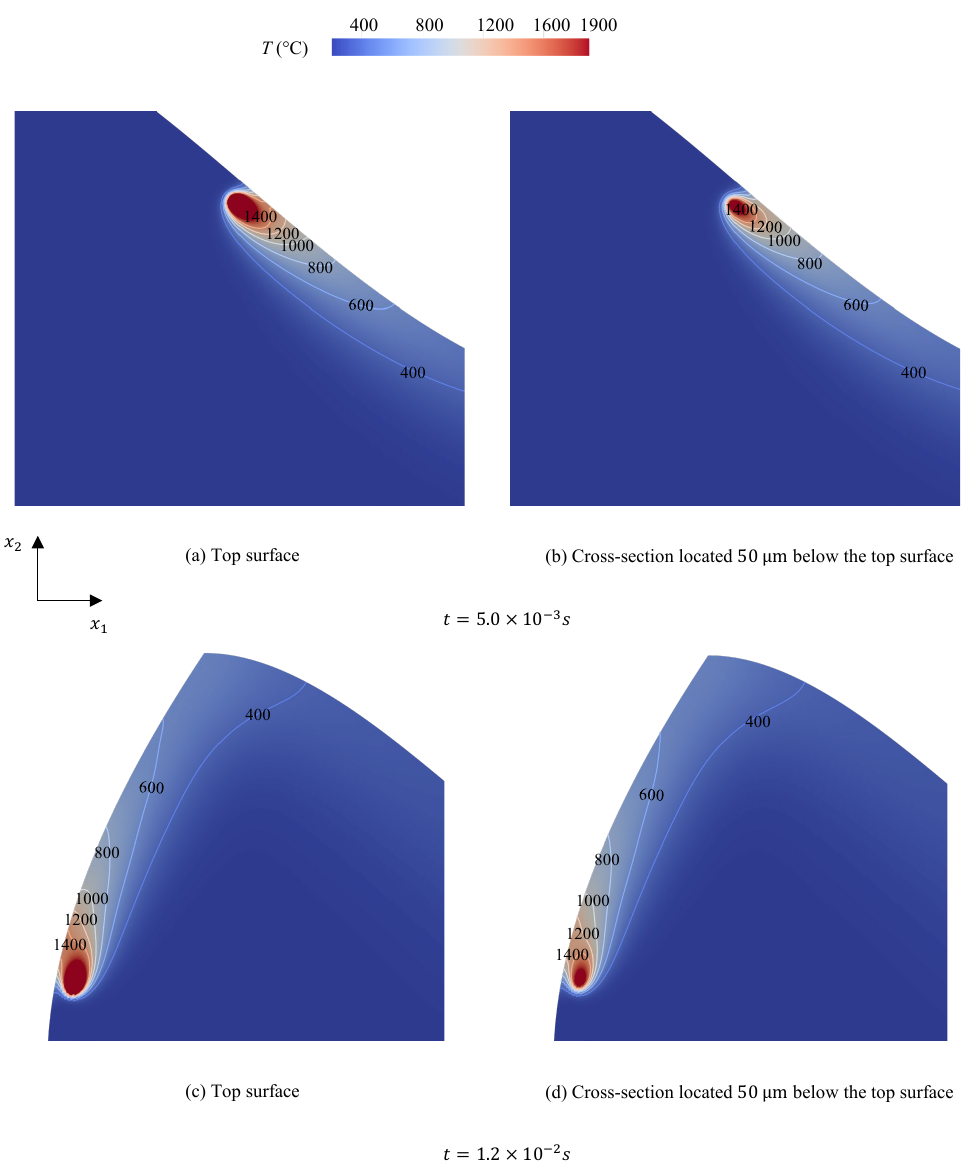}
  \caption{IGA temperature field snapshots of the top surface and a cross-section located \SI{50}{\micro\meter} below the top surface for the butterfly-shaped part at $t=\SI{5.0e-3}{s}$ and $t=\SI{1.2e-2}{s}$ during the laser contour scan.}
  \label{fig:butterfly-temperature-crosssection}
\end{figure}

Temperature distributions at $t=\SI{5.0e-3}{s}$ and $t=\SI{1.2e-2}{s}$ are presented in \cref{fig:butterfly-temperature}. At both time steps, the temperature iso-surface remain orthogonal to the part boundary surface, confirming that the adiabatic boundary conditions are properly enforced. 
The corresponding temperature contours of the top surface and of a cross section located \SI{50}{\micro\meter} below the top surface are shown in \cref{fig:butterfly-temperature-crosssection}. On both surfaces, the temperature contour lines exhibit an orientation orthogonal to the part boundaries. It is also evident that the cross-section located \SI{50}{\micro\meter} below the top surface contains a smaller region of high temperature.

\section{Conclusions and outlook}
\label{sc:conclusions}

This study has introduced and validated a novel semi-analytical framework for thermal simulation of laser powder bed fusion, which synergistically combines an analytical point-source solution with a complementary numerical field computed using isogeometric analysis. By superposing the analytical temperature field of discrete laser exposures with the IGA-based boundary correction, the method effectively enforces adiabatic boundary conditions without the need for locally refined or scan-wise adaptive meshing.

Numerical examples demonstrate the following key advantages of the proposed framework:
\begin{itemize}
    \item \textbf{Superior computational efficiency}: For a point source near a curved boundary, IGA achieved accuracy comparable to a highly refined FEM reference solution while utilizing orders of magnitude fewer degrees of freedom. Crucially, IGA maintained errors well below typical engineering tolerances even with an element size five times the laser spot radius, whereas FEM required element sizes smaller than half the spot radius for similar accuracy.
    \item \textbf{Robust boundary enforcement}: During continuous contour scanning along both simple and complex butterfly-shaped geometries, the IGA-based method reproduced reference temperature distributions with significantly coarser meshes. The consistent orthogonality of temperature iso-lines to the part boundaries throughout the scanning process confirms the method's reliability in enforcing adiabatic conditions, even for challenging non-extruded geometries with nonlinear cross-sectional variations.
    \item \textbf{Intrinsic geometric flexibility}: The spline-based IGA discretization naturally conforms to complex boundaries without resorting to impractical image-source reflections or boundary-fitted remeshing strategies, which are particularly problematic for geometries featuring sharp corners and disconnected cross-sections.
\end{itemize}

Collectively, these findings establish that the integration of semi-analytical source superposition with IGA provides a robust and computationally efficient alternative to conventional approaches. The proposed framework eliminates the dependency on dynamic mesh refinement that tracks the moving heat source while dramatically reducing sensitivity to geometric complexity, thereby offering a scalable pathway for part-scale thermal modeling of realistic LPBF components.

Future research will focus on several promising extensions: (1) enhancing the physical model by coupling the thermal field with phase-change and melt-pool fluid dynamics; (2) improving computational efficiency through adaptive local refinement of the IGA basis or hierarchical spline techniques for large-scale parts and multi-laser systems; (3) integrating the framework with process parameter optimization and feedback control for real-time thermal management; and (4) validating the model against \textit{in-situ} temperature measurements across a wider range of materials and geometries.

%%
%% The acknowledgments section is defined using the "acks" environment
%% (and NOT an unnumbered section). This ensures the proper
%% identification of the section in the article metadata, and the
%% consistent spelling of the heading.
\section*{Acknowledgments}
The authors would like to thank Prof.~Fred van Keulen for his valuable suggestions. Yang Yang is supported by China Scholarship Council (No.~202106150031).

%% bibliography
\bibliographystyle{elsarticle-num}
\bibliography{cas-refs}

\end{document}